\documentclass[MSNbibl,number,citesort,seceqn,dvips]{arxbj}
\usepackage{upgreek}

\aid{0}
\volume{19}
\issue{5B}
\pubyear{2013}
\firstpage{2389}
\lastpage{2413}
\doi{10.3150/12-BEJ456} 

\makeatletter

\newtheorem{cor}{Corollary}[section]
\newremark{ex}{Exampfle}[section]
\newtheorem{lem}{Lemma}[section]
\newtheorem{prop}{Proposition}[section]
\newremark{rem}{Remark}[section]
\newtheorem{thmm}{Theorem}[section]
\newproclaim{ass}{Assumptions}

\newcommand{\cC}{\mathcal{C}}
\newcommand{\cD}{\mathcal{D}}
\newcommand{\cG}{\mathcal{G}}
\newcommand{\cI}{\mathcal{I}}
\newcommand{\cN}{\mathcal{N}}
\newcommand{\cQ}{\mathcal{Q}}

\def\b{\beta}
\def\d{\dot}
\def\del{\delta}
\def\ep{\varepsilon}
\def\iny{\infty}

\def\si{\sigma}

\def\R{\mathbb{R}}
\def\Z{\mathbb{Z}}

\def\i{\mathrm{i}}
\makeatother

\begin{document}
\begin{frontmatter}

\title{On asymptotic distributions of
weighted sums of periodograms}
\runtitle{CLT for quadratic forms}

\begin{aug}
\author[1]{\fnms{Liudas} \snm{Giraitis}\thanksref{1}\ead[label=e1]{l.giraitis@gmail.com}}%
\and
\author[2]{\fnms{Hira L.} \snm{Koul}\corref{}\thanksref{2}\ead[label=e2]{koul@stt.msu.edu}}
\runauthor{L. Giraitis and H.L. Koul} 
\address[1]{School of Economics and Finance, Queen Mary, University of London,
Mile End Rd., London, E1 4NS, UK. \printead{e1}}
\address[2]{Department of Statistics and Probability,
Michigan State University, East Lansing, MI 48824, USA. \printead{e2}}
\end{aug}

\received{\smonth{7} \syear{2011}}
\revised{\smonth{5} \syear{2012}}

%
\begin{abstract}
We establish asymptotic normality of weighted sums of
periodograms of a stationary linear process where weights depend on the sample
size. Such sums appear in numerous statistical applications and can be
regarded as
a discretized versions of quadratic forms involving integrals of weighted
periodograms. Conditions for asymptotic normality of these weighted sums
are simple, minimal, and resemble Lindeberg--Feller
condition for weighted sums of independent and identically distributed
random variables.
Our results are applicable to
a large class of short, long or negative memory processes. The proof is
based on sharp bounds derived for
Bartlett type approximation of these sums by the corresponding sums of
weighted periodograms of independent and identically distributed random
variables.
\end{abstract}

%
\begin{keyword}
\kwd{Bartlett approximation}
\kwd{Lindeberg--Feller}
\kwd{linear process}
\kwd{quadratic forms}
\end{keyword}

\end{frontmatter}

\section{Introduction}\label{sec1}
Let $X_j, j=0,\pm1,\ldots,$ be a stationary process with a spectral
density $f_X$
and let $u_j=2\uppi j/n$, $j=1,\ldots,[n/2]$, denote discrete Fourier
frequencies.
In this paper, we develop asymptotic distribution theory for the
weighted sums
%
\begin{equation}
\label{esum11} Q_{n,X}:= \sum_{j=1}^\nu
b_{n,j}I_X(u_j), \qquad \nu:=[n/2]-1,\qquad n\ge1,
\end{equation}
of periodograms $I_X(u_j)=(2\uppi n)^{-1}  |\sum_{t=1}^n
\mathrm{e}^{\i tu_j}X_t |^2$, where $b_{n,j}$ are triangular arrays of real weights.
When $b_{n,j}=b_n(u_j),$ where
$b_n,
n\ge1$ is a sequence of real valued functions on
$\Pi:=[-\uppi,\uppi]$,
$Q_{n,X}$ is an estimate of
$\sum_{j=1}^\nu b_n(u_j)f_X(u_j)$ and
can be
viewed as a discretized version of the integral
$ 
\cI_n:= \int_0^\uppi b_n(u) I_X(u) \,\mathrm{d}u.
$ 
Integrals $\cI_n$ arise
naturally in many situations in statistical inference. For example,
the auto-covariance
function of~$\{X_j\}$ is
\[
\operatorname{ Cov}(X_k,X_0)=2\int
_0^\uppi \cos(ku)f_X(u) \,\mathrm{d}u,\qquad k=0, 1,2,
\ldots,
\]
and the spectral distribution function can be written as
$F(y)=\int_{-\uppi}^\uppi I(u\le y)f_X(u) \,\mathrm{d}u$. In these two
examples $b$ does not depend on $n$. If one wishes to estimate
$f_X(u_0)$ at a point $u_0\in(0,\uppi)$ by kernel smoothing
method,
then $b$ will typically depend on $n$.

Asymptotic distribution theory of $\cI_n$ when $b$ does not depend on
$n$ and $\{X_j\}$ is a stationary Gaussian or linear
process
is well understood and investigated both for short
memory and long memory linear processes; for asymptotic normality results
see Hannan \cite{Han73}, Fox and
Taqqu \cite{FoxTaq87}, Giraitis and Surgailis \cite{GirSur90} and Giraitis and Taqqu
\cite{GirTaq98}; for non-Gaussian limits
see
Terrin and Taqqu \cite{TerTaq90}
and Giraitis, Taqqu and Terrin \cite{GirTaqTer98}. Nualart and Peccati \cite{NuaPec05} give
simple sufficient
conditions for central limit theorem (CLT) of quadratic forms
that can be written as a sequence of multiple stochastic integrals.

It
is perhaps worth pointing out that even in the case when
$b$ does not depend on $n$, investigation of limit distribution of $\cI_n$
is technically involved.
As is evident from the works of Hannan~\cite{Han73} and
Bhansali, Giraitis and Kokoszka \cite{BhaGirKok07N2}, deriving asymptotic
distribution of $\cI_n$ 
in case of general weight sequences $b_n$ depending on $n$
will be
prohibitively complicated, and conditions for asymptotic
normality will lack desirable simplicity.

In comparison, the verification of asymptotic normality of
weighted sums of periodograms is relatively simple.
In Sections \ref{speriodo} and \ref{sgeneral} below, we provide theoretical
tools to establish the CLT for
$Q_{n,X}-EQ_{n,X}$ and $\cD_n:=
Q_{n,X}- 
\sum_{j=1}^{[n/2]}b_{n,j}f_X(u_j),
$ 
and to evaluate the large sample behavior of $E\cD_n$, $\operatorname
{Var}(Q_{n,X})$ and the
mean-squared error $E\cD_n^2$, when $\{X_j\}$ is a stationary linear process
with i.i.d. innovations, possibly having
long memory.
Our conditions for asymptotic normality of these weighted sums
are formulated in terms of $\{b_{n,j}, f_X(u_j)\}$.
They are simple and resemble Lindeberg--Feller
type condition for weighted sums of i.i.d. r.v.'s, regardless of
whether 
$\{X_j\}$ has 
short, long or negative memory.


A number of papers in the literature deal with more
general quadratic forms (sums of weighted periodograms).
Generalizations 
usually include relaxing
assumption of linearity of $\{X_j\}$. Hsing and Wu \cite{HsiWu04} obtain
asymptotic normality of a quadratic form $\sum_{ t,s=1}^n
b_{t-s}K(X_t,X_s)$ for a non-linear transform $K$ of a linear
process $\{X_j\}$ under a set of complex conditions that do not provide
a direct
answer in terms of $\{b_t\}$, $K$ and $\{X_j\}$. 
Moreover, their weights $b_t$'s are not allowed to depend on $n$.
Wu and Shao \cite{WuSha07} derive CLT for discrete Fourier transforms
and spectral density estimates under some restrictions on
dependence structure of $\{X_j\}$
based on conditional moments. Shao and Wu \cite{ShaWu07N1}
establish the CLT for quadratic forms with weights depending on
$n$ using martingale approximation method. Liu and Wu \cite{LiuWu10} consider
non-parametric
estimation of spectral density of a stationary process using
$m$-dependent approximation of $X_j$'s.
Generality of these papers requires
verification of a number of complex technical conditions
which impose {a priori} a rate condition in approximations,
that must
be verified in each specific case. For example, Wu and Shao
\cite{WuSha07} requires geometric-contraction condition, which implies
exponential decay of the autocovariance function $\gamma_X(k)$ of
$\{X_j\}$, whereas in Liu and Wu \cite{LiuWu10} the dependence is restricted
assuming summability of $|\gamma_X(k)|$.
Both papers also restrict the set of $b_{n,j}$'s to
specific weights appearing in kernel estimation. Such structural
assumptions may be easier to verify than verifying mixing conditions, but
they are redundant, not informative and too restrictive in the
case when $\{X_j\}$ is a linear process.\looseness=-1

The present paper establishes the CLT
for $Q_{n,X}$ in the latter case under minimal conditions,
which allow for short, long or negative memory 
in $\{X_j\}$ and arbitrary
weights $b_{n,j}$ as along as $f_X(u_j)b_{n,j}$'s satisfy condition
(\ref{eLF-1})
of uniform negligibility.
The main tool of the proof is
Bartlett type approximation for discrete Fourier
transforms of $X_j$'s which is essentially different from
the methods of approximations used in the above works.
Besides being simple and easy to verify, the obtained conditions
are close to being necessary; see Remark \ref{rtq} below.


\begin{ass*} Accordingly, let $\Z:=\{0,\pm
1,\ldots\}$,
%
\begin{equation}
\label{elinear} X_j=\sum_{k=0}^{\infty}a_k
\zeta_{j-k}, \qquad j\in\Z,\ \sum_{k=0}^{\infty}a_k^{2}<
\infty,
\end{equation}
be a
linear process where $\{\zeta_j, j\in\Z\}$ are i.i.d. standardized r.v.'s.
Assume that the spectral density
$f_X$ of the process $X_j, j\in\Z,$
satisfies
%
\begin{eqnarray}
f_X(u )=|u |^{-2d}g(u ), \qquad |u |\leq\uppi, \label{effNew}
\end{eqnarray}
for some
$|d|<1/2$, where $g(u)$ is a continuous function satisfying
\[
0<C_1\le g(u)\le C_2<\infty,\qquad u\in\Pi\ (\exists
0<C_1,C_2<\infty).
\]
\end{ass*}

Condition (\ref{effNew}) allows to derive the mean square error
bounds of estimates, which are given in Theorem
\ref{thmbias-t}. To derive asymptotic normality and some delicate
Bartlett type approximations, we shall additionally need to assume that
the transfer function
$A_X(u):=\sum_{k=0}^\infty \mathrm{e}^{-\i k u}a_k$, $u\in\Pi$,
is differentiable in $(0,\uppi)$ and its derivative $\d A_X$ satisfies
%
\begin{equation}
\bigl|\dot{A}_X(u) \bigr|\le C|u|^{-1-d},\qquad  u \in\Pi. \label{eaaa}
\end{equation}


Conditions (\ref{effNew}) and (\ref{eaaa}) are formulated this
way to cover long and negative memory models, with $|d|<1/2$,
$d\ne0$. They allow spectral density to vanish or to have a
singularity
point at zero frequency. The short memory case where
$f_X$ and $A_X$ are Lipshitz continuous and bounded away from $0$ and
$\infty$
is also discussed in Section \ref{sgeneral}.

To proceed further,
define the \textit{discrete Fourier transforms} (DFT) of $\{X_j\}$ and
$\{\zeta_j\}$ computed at frequencies $u_j$'s,
$j=0,\ldots,[n/2],$ to be, respectively,
\[
 w_{X,j}=\frac{1}{\sqrt{2\uppi
n}}\sum
_{k=1}^n \mathrm{e}^{\i u_j k}X_k,\qquad
w_{\zeta,
j}=\frac{1}{\sqrt{2\uppi n}}\sum_{k=1}^n
\mathrm{e}^{\i u_j k}\zeta_k.
\]
%
The corresponding periodograms, transfer
functions and spectral densities of $\{X_j\}$ and $\{\zeta_j\}$
at frequency $u_j$
are denoted by
\begin{eqnarray*}
I_{X,j}&=&|w_{X,j}|^2,\qquad I_{\zeta,j}=|w_{\zeta, j}|^2,\qquad
A_{X,j}=A_{X}(u_j), \qquad A_{\zeta,j}=1,
\\
f_{X,j}&:=&f_X(u_j),\qquad f_{\zeta,j}:=f_\zeta(u_j)
\equiv\frac
{1}{2\uppi}, \qquad j=0,1,\ldots,[n/2]. \label{eperod}
\end{eqnarray*}

The goal of establishing asymptotic normality of
$Q_{n,X}$ is facilitated by first developing
asymptotic distribution theory for the sums
\[
S_{n,X}:=\sum_{j=1}^{\nu}
b_{n,j}\frac{I_{X,j}}{f_{X,j}}.
\]
Moreover, asymptotic analysis of these sums is more illustrative of the
methodology used. The asymptotic normality of
$S_{n,X}$ is discussed in Section \ref{speriodo}.

The CLT for the quadratic forms $Q_{n,X}$
with weights 
not depending on $n$
was investigated by Hannan \cite{Han73}; see also
Proposition 10.8.6. of Brockwell and Davis \cite{BroDav91}. Their proof
required restrictive condition $\sum_{k=0}^\infty
k^{1/2}|a_k|<\infty$ on the coefficients $a_k$ of the linear process
$\{X_j\}$ of (\ref{elinear})
and was based on Bartlett approximation of periodogram
$I_{X,j}/f_{X,j}$ by periodogram $I_{\zeta,j}/f_{\zeta,j}$ of the
noise. The idea for
the theory and the proofs presented in this paper have their roots in
Robinson \cite{Rob95N2}.

We show that CLT's for $Q_{n,X}$ and $S_{n,X}$
hold under similar conditions as the classical CLT for weighted sums of
i.i.d. r.v.'s.
It requires Lindeberg--Feller type condition on weights
$b_{n,j}$ and minimal restrictions on a linear process $\{X_j\}$
which may have short or long memory. For example,
in short memory case it suffices to assume that
$a_k$ of 
(\ref{elinear}) satisfy
$\sum_{k=0}^\infty|a_k|<\infty$ and 
$f_X$ is bounded away from $0$ and $\infty$; see Section \ref{sgeneral}.
Results below also show that weighted sums of rescaled periodogram
$I_{X,j}/f_{X,j}$ of a linear process behave, to some extend, similarly
as the
weighted sums of i.i.d. r.v.'s.

We also investigate precision of Bartlett
approximation of $Q_{n,X}$ and $S_{n,X}$ by sums of
weighted periodograms
$I_{\zeta,j}/f_{\zeta,j}$. 
Lemma \ref{lemlap} and Theorem \ref{thmbias-t}
contain sharp bounds and are of independent interest. From these results,
one sees that the above approximation is extremely precise, and the
resulting error is small and can be effectively controlled by
the weights $\{b_{n,j}\}$ alone.
This type of approximation is a popular tool for
establishing CLT for specific types of weights
$b_{n,j}$, for example,
for local Whittle estimators; see Robinson
\cite{Rob95N2},
Shao and Wu \cite{ShaWu07N2} and Shao~\cite{Sha10}.
In these papers, innovation sequence is allowed to be
a martingale difference
or an uncorrelated weakly dependent non-linear causal process.
However, because of narrower focus, they deal with special
weights and do not seek establishing a general CLT for
$Q_{n,X}$ as such. 
In our setting, assumption of i.i.d. innovations is a secondary
issue and also can be relaxed, while the major objective is
obtaining the CLT for $Q_{n,X}$ with the most general feasible
weighting scheme $b_{n,j}$.

Finally, in the present paper the spectral density $f_X$ is
allowed to take infinite or zero value only at the zero frequency
restricting $|d|<1/2$ to keep $\{X_j\}$ stationary. Establishing
sufficient conditions for CLT
for a differenced stationary process,
as well as when the spectral density $f_X$ may have
singularity/zero at a frequencies away from zero is of definite
interest, but needs further investigation.

In the sequel, Cum$_k(Z)$ denotes the
$k$th cummulant of the r.v. $Z$, $\operatorname{IID}(0,1)$ denotes the class of i.i.d.
standardized
r.v.'s, $a\wedge b:=\min(a,b)$, $a\vee b:=\max(a,b)$, for any real
numbers $a, b$, and
all limits are taken s $n\to\iny$, unless specified otherwise.


\section{Asymptotic normality 
of $S_{n,X}$}\label{speriodo}

Important role in the asymptotic analysis
of $S_{n,X}$ is played by Bartlett type approximation
\[
(I_{X,j}/f_{X,j} )\sim (I_{\zeta,j}/f_{\zeta,j} )=2
\uppi I_{\zeta,j},\qquad  j=1,\ldots,\nu, \nu=[n/2]-1.
\]
Our first goal is to
approximate $S_{n,X}$ by the weighted sum of 
$I_{\zeta,j}$,
%
\begin{equation}
\label{esnz}S_{n,\zeta} =\sum_{j=1}^{\nu}
b_{n,j} (I_{\zeta,j}/f_{\zeta,j} ) \equiv \sum
_{j=1}^{\nu} b_{n,j} 2\uppi I_{\zeta,j}.
\end{equation}
Let
%
\begin{eqnarray}
\label{edef-a} R_n&:=&S_{n,X}-S_{n,\zeta},\qquad
b_n := \max_{j=1,\ldots, {\nu}} |b_{n,j}|,\qquad B_n:=
\Biggl(\sum_{j=1}^{\nu} b_{n,j}^2
\Biggr)^{1/2},
\nonumber
\\[-8pt]
\\[-8pt]
\nonumber
q_n^2&:=&B_n^2+\operatorname{Cum}_4(\zeta_0)\frac{1}{n} \Biggl(\sum
_{j=1}^{\nu
}b_{n,j} \Biggr)^2.
\end{eqnarray}
We show later that $\operatorname {Var}(S_{n,\zeta})=q_n^2$; see (\ref{es-mean})(b).

Lemma \ref{lemlap} below
provides an upper bound of order $b_n \log^2(n)$ for $ER_n^2$ while
Lemma \ref{lemlasymp}
establishes the asymptotic normality of the approximating quadratic
forms $S_{n,\zeta}$.
The following theorem gives an approximation to $ES_{n,X}$, $\operatorname
{Var}(S_{n,X})$, and
proves asymptotic normality of $S_{n,X}$
under Lindeberg--Feller type condition (\ref{es-con1}) on the weights
$b_{n,j}$.

Because of the invariance property 
$I_{X+\mu}(u_j)=I_X(u_j)$, $\mu\in\R$, $j=1,\ldots, n-1$,
all results obtained below
remain valid also for a process $\{X_j\}$ of (\ref{elinear}) that has
non-zero mean.
%
\begin{thmm}\label{ts-main} Suppose the linear process $\{X_j, j\in
\Z\}$
of (\ref{elinear}) satisfies assumptions
(\ref{effNew}) and~(\ref{eaaa}), and $E\zeta^4_0$ $<\infty$.
About the weights $b_{n,j}$'s
assume
%
\begin{equation}
\frac{b_n}{B_n}=\frac{\max_{j=1,\ldots, \nu}|b_{n,j}|}{ (\sum_{j=1}^{\nu}b^2_{n,j} )^{1/2}} \to0. \label{es-con1}
\end{equation}
Then
%
\begin{eqnarray}
\label{es-c}  E S_{n,X} &=& \sum_{j=1}^{\nu}
b_{n,j}+\mathrm{o}(q_n), \qquad\operatorname{ Var}(S_{n,X})=
q_n^2+\mathrm{o}\bigl(q_n^2\bigr),
\nonumber
\\[-8pt]
\\[-8pt]
\nonumber
\operatorname{ Var}(S_{n,X})^{-1/2} (S_{n,X}-ES_{n,X}
) &\to_D &\cN(0,1),\qquad q_n^{-1}
\Biggl(S_{n,X}-\sum_{j=1}^{\nu}
b_{n,j} \Biggr) \to_D \cN(0,1).
\end{eqnarray}
Moreover,
%
\begin{equation}
\label{ebab14}\min \bigl(1,\operatorname{Var}\bigl(\zeta_0^2
\bigr)/2 \bigr)B_n^2\le q_n^2\le
\bigl(1+\bigl|\operatorname{ Cum}_4(\zeta_0)\bigr|
\bigr)B_n^2.
\end{equation}
\end{thmm}
\begin{pf} The proof uses Lemmas \ref{lemlap} and \ref{lemlasymp}
given below. To prove (\ref{ebab14}), use definition of
$q_n$ and the Cauchy--Schwarz inequality to obtain the
upper bound. The lower bound is derived in (\ref{ecum4}) of Lemma
\ref{lemlasymp}.

By (\ref{es-con1}), (\ref{ebab14}), (\ref{esacon})(b), and (\ref
{es-mean}),
%
\begin{equation}
\label{eraw} ES_{n,\zeta}=\sum_{j=1}^{\nu}
b_{n,j},\qquad  E|R_n|\le\bigl(ER^2_n
\bigr)^{1/2}=\mathrm{o}(B_n)=\mathrm{o}(q_n).
\end{equation}
These facts 
in turn complete the proof of the first claim in (\ref{es-c}).

To prove the second claim,
note that by (\ref{es-mean})(b), $\operatorname{ Var}(S_{n,\zeta})=q_n^2$, which
together with (\ref{eraw}) yields
$\operatorname{ Var}(R_n)\le ER_n^2
=\mathrm{o}(q_n^2),$ 
$|\operatorname{ Cov}(S_{n,\zeta},R_n)|
=\mathrm{o}(q_n^2).$
These facts together with a routine argument
complete the proof of the second claim in (\ref{es-c}).

Finally, again in view of (\ref{eraw}),
%
\[
S_{n,X}-\sum_{j=1}^{\nu}
b_{n,j}= S_{n,X}-ES_{n,\zeta} 
=S_{n,\zeta}-ES_{n,\zeta}+\mathrm{o}_p(q_n).
\]
This and (\ref{es-mean})(c) of Lemma \ref{lemlasymp} imply the
first asymptotic
normality result in (\ref{es-c}), while the last claim follows from
the first three claims in
(\ref{es-c}).
\end{pf}

Lemma \ref{lemlap} below provides the two types of sharp upper
bounds for
$ER_n^2$ that are useful in approximating
$S_{n,X}$ by $S_{n,\zeta}$.
The idea of using Bartlett type approximations to establish the
asymptotic normality of an integrated weighted periodogram of a
short memory linear process
goes back to the work of Grenander and Rosenblatt \cite{GR53}, Hannan and Heyde
\cite{HanHey72} and Hannan~\cite{Han73}, whereas for sums of weighted
periodograms of an ARMA process it was used in Proposition~10.8.5
of Brockwell and Davis~\cite{BroDav91}. Their approximations were derived
under the assumption that the weight function $b$ did not depend
on $n$, and the bounds they obtain
have low-level of sharpness,
though they are sufficient to show that the main term dominates the remainder.
The sharp bounds
for an integrated weighted periodogram established in Bhansali
\textit{et al}.~\cite{BhaGirKok07N2} technically are more involved
and harder to apply than those for sums in this lemma.
%
\begin{lem}\label{lemlap} Assume that $\{X_j\}$ of (\ref{elinear})
satisfies
(\ref{effNew}) and (\ref{eaaa}), and $E\zeta^4_0<\infty$. Then
%
\begin{eqnarray}
\label{es-con0>} E(R_n-ER_n)^2&\le& C
b_n^2 \log^{3}(n)\quad \mbox{and}\quad
E(R_n-ER_n)^2 \le Cb_n
B_n,
\\
\label{es-con0?} |ER_n|&\le& C b_n \log^2(n)\quad
\mbox{and} \quad|ER_n| = \mathrm{o}(B_n)\qquad \mbox{if } b_n
=\mathrm{o}(B_n).
\end{eqnarray}
In particular,
%
\begin{eqnarray}
\label{esacon} \mathrm{(a)}\quad E(S_{n,X}-S_{n,\zeta})^2&
\le& C b^2_n \log^{4}(n);
\nonumber
\\[-8pt]
\\[-8pt]
\nonumber
\mathrm{(b)}\quad  E(S_{n,X}-S_{n,\zeta})^2 &= & o
\bigl(B_n^2\bigr) \qquad\mbox{if } b_n
=\mathrm{o}(B_n).
\end{eqnarray}
\end{lem}

The proof of this lemma is facilitated by the following two propositions.
%
\begin{prop}
\label{pbound-f} Let $\{Y_{n,j}^{(i)}, j=1,\ldots, n\}$,
$i=1,2$, $n \ge1$
be the two sets of moving averages
\begin{eqnarray*}
Y_{n,j}^{(i)}=\sum_{k\in\Z}
b^{(i)}_{n,j}(k)\zeta_{k}, \qquad\sum
_{k\in\Z}^\infty\bigl|b^{(i)}_{n,j}(k)\bigr|^2<
\infty,\qquad i=1,2,
\end{eqnarray*}
where $\{b^{(i)}_{n,j}(k)\}$ are possibly complex weights.
Assume, $\zeta_k\sim \operatorname{IID}(0,1)$,
$E\zeta_0^4<\infty$.
Then, for any real weights $c_{n,j},$ $ j=1,\ldots, n$,
%
\begin{eqnarray}
\label{eriba-va1} &&
 \operatorname{ Var} \Biggl(\sum
_{j=1}^{ n} c_{n,j}\bigl\{\bigl|Y_{n,j}^{(1)}\bigr|^2-\bigl|Y_{n,j}^{(2)}\bigr|^2
\bigr\} \Biggr)
\nonumber
\\[-8pt]
\\[-8pt]
\nonumber
&&\quad\le \bigl(4+4\operatorname{ Var}\bigl(\zeta_0^2\bigr)
\bigr) 
\sum_{j,k=1}^{ n}
|c_{n,j}c_{n,k}| \bigl| \bigl|r^{11}_{n,jk}\bigr|^2+\bigl|r^{22}_{n,jk}\bigr|^2
-2\bigl|r^{12}_{n,jk}\bigr|^2 \bigr|,
\end{eqnarray}
where $
r^{il}_{n,jk}:=E[Y_{n,j}^{(i)} \overline{Y_{n,k}^{(l)}}]=
\sum_{t\in\Z}b^{(i)}_{n,j}(t) \overline{b^{(l)}_{n,k}(t)},\
i,l=1,2. $
\end{prop}
\begin{pf} Observe that
\begin{eqnarray*}
\cG_n &:=& \sum_{j=1}^{ n}
c_{n,j}\bigl\{ \bigl|Y_{n,j}^{(1)}\bigr|^2-\bigl|Y_{n,j}^{(2)}\bigr|^2
\bigr\}
\\
&=&\sum_{t, s\in\Z} \Biggl(\sum
_{j=1}^nc_{n,j}\bigl\{b^{(1)}_{n,j}(t)
\overline{b^{(1)}_{n,j}(s)}-b^{(2)}_{n,j}(t)
\overline {b^{(2)}_{n,j}(s)}\bigr\} \Biggr)
\zeta_{t}\zeta_{s} =: \sum_{t, s\in\Z}B_n(t,s)
\zeta_{t}\zeta_{s}.
\end{eqnarray*}
Hence,
\begin{eqnarray*}
&&E |\cG_n-E\cG_n |^2\\ 
&&\quad\le 4 \biggl(E \biggl|\sum_{t< s }B_n(t,s)
\zeta_{t}\zeta_{s} \biggr|^2 +E \biggl|\sum
_{ s<t}B_n(t,s)\zeta_{t}
\zeta_{s} \biggr|^2 +E \biggl|\sum_{t\in\Z} B_n(t,t) \bigl(
\zeta_{t}^2-E\zeta_t^2\bigr)
\biggr|^2 \biggr)
\\
&&\quad= 4\sum_{ t< s}\bigl|B_n(t,s)\bigr|^2
+4\sum_{ s<t}\bigl|B_n(t,s)\bigr|^2
+4\operatorname{ Var}\bigl(\zeta_0^2\bigr)
\sum_{t\in\Z}\bigl|B_n(t,t)\bigr|^2
\\
&&\qquad\le\bigl(4+4\operatorname{ Var}\bigl(\zeta_0^2\bigr)
\bigr)\sum_{t, s\in\Z}\bigl|B_n(t,s)\bigr|^2.
\end{eqnarray*}
But,
\begin{eqnarray*}
&&\sum_{t, s\in\Z}\bigl|B_n(t,s)\bigr|^2\\
&&\qquad=
\sum_{j,k=1}^nc_{n,j}c_{n,k}
\sum_{t, s\in\Z}\bigl\{b^{(1)}_{n,j}(t)
\overline{b^{(1)}_{n,j}(s)}-b^{(2)}_{n,j}(t)
\overline {b^{(2)}_{n,j}(s)}\bigr\}\bigl\{\overline{b^{(1)}_{n,k}(t)}
b^{(1)}_{n,k}(s)-\overline{b^{(2)}_{n,k}(t)}
b^{(2)}_{n,k}(s)\bigr\}
\\
&&\qquad=\sum_{j,k=1}^nc_{n,j}c_{n,k}
\bigl(\bigl|r^{11}_{n,jk}\bigr|^2+\bigl|r^{22}_{n,jk}\bigr|^2-\bigl|r^{12}_{n,jk}\bigr|^2
-\bigl|r^{12}_{n,kj}\bigr|^2 \bigr).
\end{eqnarray*}
This completes the proof of (\ref{eriba-va1}).
\end{pf}

The next proposition describes some needed large sample properties
of DFTs. Because $ \sum_{t=1}^ne^{\i tu_m}
=n\{I(m=0)+I(m=n)\},$ DFTs of a white noise process $\{\zeta_j\}$
are uncorrelated:
%
\begin{eqnarray}
\label{epro-dtf} E[w_{\zeta,j}\overline{w_{\zeta,k}}]&=&
\frac{1}{2\uppi},\qquad 1\le k=j\le n,
\nonumber
\\[-8pt]
\\[-8pt]
\nonumber
&=& 0, \qquad 1\le k<j\le n.
\nonumber
\end{eqnarray}
%

Consider now the two linear processes
\[
X_j=\sum_{k=0}^{\infty
}a_k
\zeta_{j-k},\qquad Y_j=\sum_{k=0}^{\infty}b_k
\zeta_{j-k},\qquad j\in\Z,\qquad \sum_{k=0}^{\infty}a_k^{2}<
\infty,\qquad \sum_{k=0}^{\infty}b_k^{2}<
\infty,
\]
with the same white noise
innovations $\{\zeta_j\}\sim WN(0,\si^2)$.
Let $ 
A_X(v):=\sum_{k=0}^\infty \mathrm{e}^{-\i kv}a_k,$ %
$A_Y(v):=\sum_{k=0}^\infty \mathrm{e}^{-\i kv}b_k,
$ 
$f_X(v)=(\sigma^2/2\uppi)|A_X(v)|^2,$ %
$f_Y(v)=
(\sigma^2/2\uppi)|A_Y(v)|^2, $
denote their respective transfer and 
spectral densities.

Let $f_{XY}(v)$ denote a (complex valued) cross-spectral
density:
%
\begin{eqnarray}
\label{ecross-sp} f_{XY}(v)&:=&\frac{\sigma^2}{2\uppi} A_X(v)
\overline{A_Y(v)},\qquad v\in\Pi,
\nonumber
\\[-8pt]
\\[-8pt]
\nonumber
E[X_jY_{j-k}]&=&\int_\Pi
\mathrm{e}^{\i kv}f_{XY}(v)\,\mathrm{d}v =\frac{\sigma^2}{2\uppi}\sum
_{l=0}^{\infty}a_{l+k}b_{l},\qquad  k\ge0,
j\in\Z.
\end{eqnarray}
If $Y_j=\zeta_j, $ $j\in\Z$, then
\begin{eqnarray*}
\label{ecross-Af} f_{X\zeta}(v)&:=&\frac{\sigma^2}{2\uppi} A_X(v),\qquad
v\in\Pi, \\
E[X_j \zeta_{j-k}]&=&\frac{\sigma^2}{2\uppi}\int
_\Pi \mathrm{e}^{\i
kv}A_X(v)\,\mathrm{d}v=
\sigma^2a_k,\qquad  k\ge0.
\end{eqnarray*}

Proposition \ref{propf-mult} below summarizes asymptotic properties of
cross-covariances $E[w_{X,j}\overline{w_{Y,k}}]$. It
generalizes and extends Theorem 2 of
Robinson \cite{Rob95N1}
for short memory
and long memory time series, which enable derivation of the
upper bounds based on Bartlett approximation of this paper. Its proof is
technical and appears in Giraitis, Koul and Surgailis \cite{GirKou}.

In case when Fourier frequencies in covariances
$E[w_{X,j}\overline{w_{Y,k}}]$ are from an interval
$(-\Delta,\Delta)$, $\Delta<\uppi$ (a neighborhood of $0$), smoothness
conditions on $f_X, f_Y$, $A_X,A_Y$ are local, that is, they need to be
imposed on an interval $[0,a]$, $a>\Delta$. 

To proceed further, for any subset $A \subset\R$, let $\cC(A)$ 
denote complex valued functions that are continuous on $A$, 
and $\Lambda_\b(A)$ 
denote Lipschitz continuous
functions on $A$ 
with parameter $\beta\in(0,1]$. We write $h\in\mathcal{
C}_{1,\alpha}[0,a]$, $|\alpha|< 1$, $a>0$, if
\begingroup
\abovedisplayskip=7pt
\belowdisplayskip=7pt
\[
\bigl|h(u)\bigr|\le C|u|^{-\alpha},\qquad \bigl|\dot{h}(u)\bigr|\le C|u|^{-1-\alpha}\qquad \forall u
\in[0,a].
\]
Members of $\mathcal{ C}_{1,\alpha}[0,a]$
can have an infinite peak and can be non-differentiabile
at $0$, whereas $\Lambda_\beta[0,a]$ covers continuous piecewise
differentiable
functions.

Note that for any
$h\in\mathcal{ C}[0,a]$, 
$ \omega_h(\eta):= \sup_{u,v\in[0,a]: |u-v|\le
\eta}|h(u)-h(v)|\to0,$ as $\eta\to0$.
Define $
\delta_{n,\varepsilon}(h):=\omega_h( n^{-1}\log(n))+(\log
(n))^{-\varepsilon},$ $ 0<\ep<1.$
We also need to introduce 
%
\begin{eqnarray}
\label {ellnj} \ell_n(\varepsilon;k) &:=& \frac{\log(2+
k)}{(2+k)^{1-\varepsilon}}+
\frac{\log
(2+n-k)}{(2+n-k)^{1-\varepsilon}},\qquad  0\le k\le j\le n,
\nonumber\\[-1pt]
r_{n,jk}(g)&:=&0,\qquad  g\in\Lambda_1[0, a],\ \beta=1,
\nonumber
\\[-8pt]
\\[-8pt]
\nonumber
&:=&n^{-\beta} \ell_n(\b;j-k),\qquad  g\in\Lambda_\beta[0,
a],\ 0<\beta<1,
\\[-1pt]
&:=&\delta_{n,\varepsilon}(g) \ell_n(\varepsilon;j-k), \qquad g\in
\mathcal{ C}[0, a],\ \varepsilon\in(0, 1).
\nonumber
\end{eqnarray}

\begin{prop}\label{propf-mult} Let either $\Delta<a<\uppi$, or
$\Delta=a=\uppi$.
Then, the following facts \textup{(i)--(iv)} hold for all $0<|u_k|\le u_j<
\Delta$,
\begin{longlist}[(iii)] 
\item[(i)] If
$f_{XY}\in\Lambda_\beta[0,a]$, $0<\beta\le1$, then
\begin{eqnarray*}
\bigl|E[w_{X,j}\overline{w_{Y,j}}]-f_{XY}(u_j)
\bigr| &\le& C n^{-1}\log(n),\qquad \beta=1,
\\[-1pt]
&\le& C n^{-\beta}, \qquad 0<\beta<1.
\\[-1pt]
\bigl|E[w_{X,j}\overline{w_{Y,k}}] \bigr| &\le&
C n^{-1}\log(n), \qquad \beta=1,
\\[-1pt]
&\le& C n^{-\beta} \ell_n(\b;j-k),\qquad  0<\beta<1, k<j. 
\end{eqnarray*}

\item[(ii)] If $f_{XY}\in\mathcal{ C} [0,a]$, then,
$\forall\varepsilon
\in(0,1)$,
\begin{eqnarray*}
\bigl|E[w_{X,j}\overline{w_{Y,j}}]-f_{XY}(u_j)
\bigr| &\le& C \delta_{n,\varepsilon}(f_{XY}),
\\[-1pt]
\bigl|E[w_{X,j}\overline{w_{Y,k}}]\bigr| &\le& C \delta_{n,\varepsilon}(f_{XY})
\ell_n(\ep;j-k), \qquad k<j. 
\end{eqnarray*}
\item[(iii)] If $f_{XY}\in\mathcal{ C}_{1,\alpha} [0,a]$,
$|\alpha|<1$, then
\begin{eqnarray*}
\bigl|E[w_{X,j}\overline{w_{Y,j}}]-f_{XY}(u_j)
\bigr| &\le& C u_j^{-\alpha} j^{-1} \log(1+j),
\\[-1pt]
\bigl|E[w_{X,j}\overline{w_{Y,k}}] \bigr|&\le& C \bigl(|u_k|^{-\alpha}+u_j^{-\alpha}
\bigr) j^{-1} \log j,\qquad  k<j.
\end{eqnarray*}
\item[(iv)] Suppose $f_{XY}=hg$,
where $h\in\mathcal{ C}_{1,\alpha} [0,a]$, $|\alpha|<1$, and
$g\in\Lambda_\beta[0,a]\cup\mathcal{ C} [0,a] $, $0<\beta\le1$. Then
%
\begin{eqnarray*}
&& \bigl|E[w_{X,j}\overline{w_{Y,k}}]- f_{XY}(u_j)I(j=k)
\bigr|
\\[-1pt]
&&\qquad\le C \bigl(\bigl(|u_k|^{-|\alpha|}+u_j^{-|\alpha|}
\bigr) j^{-1}\log j + \bigl(|u_k|^{-|\alpha|}\wedge
u_j^{-|\alpha|} \bigr) r_{n,jk}(g) \bigr).
\end{eqnarray*}
\end{longlist}
The constant $C$ in the above \textup{(i)--(iv)} does not depend on
$k, j$ and $n$.
\end{prop}\eject
\endgroup

Parts (i) and (ii) of the above
proposition consider the case when
$f_X$ is continuous and bounded, whereas part (iv) covers
the case when $f_X$ satisfies (\ref{effNew}) with
a bounded and continuous $g$.
The case when $g$ has also bounded derivative is covered in
part (iii). Obtaining upper bounds in the above proposition
does not require the process $\{X_j\}$ to be
linear. For convenience of applications, this proposition is formulated
for a cross-spectral density of two stationary linear processes
with the same underlying white noise innovations.
This allows to express their cross spectral density
via their transfer functions 
as indicated
in (\ref{ecross-sp}). In general, the results of Proposition
\ref{propf-mult} are valid for any spectral density or
cross-spectral density that satisfies the assumed smoothness condition.

Lahiri \cite{Lah03} provides a characterization of
asymptotic independence of the DFTs in terms of the distance between
their arguments under both
short- and long-range dependence of the underlying process.
Nordman and Lahiri \cite{NorLah06} contains some relevant
results about Bartlett correction of the frequency domain empirical
likelihood ratios.

Now rewrite
%
\begin{eqnarray}
\label {esxsz} R_n&=&S_{n,X}-S_{n,\zeta}= \sum
_{j=1}^{ \nu} b_{n,j} \biggl(
\frac{I_{X,j}}{f_{X,j}} -\frac{I_{\zeta,j}}{f_{\zeta
,j}} \biggr) =\sum_{j=1}^{ \nu}
\frac{b_{n,j}}{f_{X,j}} \biggl\{I_{X,j}-f_{X,j}\frac{I_{\zeta,j}}{f_{\zeta,j}}
\biggr\}.
\end{eqnarray}
The corollary below, which follows from Proposition \ref{pbound-f},
is useful in analyzing the sums of the types appearing in (\ref{esxsz}).
Let $f_{X\zeta,j}:=f_{X\zeta}(u_j)$.
%
\begin{cor}\label{ccor-b} Suppose that $\{X_j\}$ is a linear process
as in (\ref{elinear}) and
$E\zeta_0^4<\infty$. Then, for any real weights
$c_{n,j}$, $j=1,\ldots, n$,
%
\begin{eqnarray}
\label{eri-co2} \operatorname{Var} \Biggl(\sum_{j=1}^{ \nu}
c_{n,j}\biggl\{I_{X,j}-f_{X,j}\frac
{I_{\zeta,j}}{f_{\zeta,j}}\biggr
\} \Biggr) \le C(s_{n,1}+s_{n,2}),
\end{eqnarray}
where
\begin{eqnarray*}
s_{n,1}&:= &C\sum_{j=1}^{ \nu}
c_{n,j}^2 \bigl\{\bigl(E|w_{X,j}|^2-f_{X,j}
\bigr)^2 +f_{X,j} \bigl|E|w_{X,j}|^2-f_{X,j}
\bigr|\\
 &&\hspace*{44pt}{}+ f_{X,j} \bigl|E[w_{X,j}\overline{w_{\zeta
,j}}]-f_{X\zeta
,j}
\bigr|^2 +f_{X,j}^{3/2} \bigl|E[w_{X,j}
\overline{w_{\zeta,j}}]-f_{X\zeta
,j} \bigr| \bigr\},
\nonumber
\\
s_{n,2}&:=&\sum_{1\le k<j\le\nu}|c_{n,j}c_{n,k}|
\bigl\{ \bigl|E[w_{X,j}\overline{w_{X,k}}]\bigr|^2+
f_{X,k}\bigl|E[w_{X,j}\overline{w_{\zeta,k}}]\bigr|^2
\bigr\}.
\end{eqnarray*}
\end{cor}

\begin{pf} Observe\vspace*{-1pt} that $f_{X,j}I_{\zeta,j}/f_{\zeta,j}=|A_{X,j}|^2
I_{\zeta,j}=|A_{X,j}w_{\zeta,j}|^2,$
and that $Y_{n,j}^{(1)}:=w_{X, j}$ and $Y_{n,j}^{(2)}:=A_{X,j}w_{\zeta
, j}$ 
are moving averages of $\zeta_j$'s with complex weights. Hence, by
Proposition~\ref{pbound-f},
the l.h.s. of (\ref{eri-co2}) 
is bounded above by
\begin{eqnarray*}
&&\hspace*{-4pt} C\sum_{j,k=1}^{ \nu} |c_{n,j}c_{n,k}|
\bigl|\bigl|E[w_{X,j}\overline{w_{X,k}}]\bigr|^2
+|A_{X,j}|^2|A_{X,k}|^2\bigl|E[w_{\zeta,j}
\overline{w_{\zeta,k}}]\bigr|^2 -2|A_{X,k}|^2\bigl|E[w_{X,j}
\overline{w_{\zeta,k}}]\bigr|^2\bigr |
\\
&&\quad = C \Biggl(\sum_{j=k=1}^{ n}[\cdots]+
\sum_{k\ne
j}[\cdots] \Biggr):=C\bigl(s^\prime_{n,1}+s^\prime_{n,2}
\bigr).
\end{eqnarray*}
By (\ref{epro-dtf}), $E|w_{\zeta,j}|^2=1/2\uppi$,
$E[w_{\zeta,j}\overline{w_{\zeta,k}}]=0$, for $1\le k<j\le\nu$. Recall
also that $f_{X,j}=|A_{X,j}|^2/(2\uppi).$ Therefore,
\begin{eqnarray*}
s^\prime_{n,1} &=&\sum_{j,k=1}^{ \nu}
c_{n,j}^2 \bigl|\bigl(E|w_{X,j}|^2
\bigr)^2+f_{X,j}^2 -4\uppi f_{X,j}\bigl|E[w_{X,j}
\overline{w_{\zeta,j}}]\bigr|^2 \bigr|,
\\
s_{n,2}^\prime&=&\sum_{1\le k<j\le\nu}|c_{n,j}c_{n,k}|
\bigl(\bigl|E[w_{X,j}\overline{w_{X,k}}]\bigr|^2+
f_{X,k}\bigl|E[w_{X,j}\overline{w_{\zeta,k}}]\bigr|^2
\bigr)=s_{n,2} .
\end{eqnarray*}
%
To bound $s_{n,1}^\prime$, let
$A:=(E|w_{X,j}|^2)^2-f_{X,j}^2,$ and $
B:=|E[w_{X,j}\overline{w_{\zeta,j}}]|^2-f_{X\zeta,j}.
$ 
Then use the fact that $4\uppi f_{X,j}|f_{X\zeta,j}|^2=4\uppi
f_{X,j}|A_{X,j}|^2/(2\uppi)^2=2f_{X,j}^2$ to
rewrite the term within $|\cdots|$ in $s_{n,1}^\prime$ as
\begin{eqnarray*}
&&\bigl(E|w_{X,j}|^2\bigr)^2+f_{X,j}^2
-4\uppi f_{X,j}\bigl|E[w_{X,j}\overline{w_{\zeta,j}}]\bigr|^2
\\
&&\quad=(A - 4\uppi f_{X,j}B)+\bigl(2f_{X,j}^2-4\uppi
f_{X,j}|f_{X\zeta,j}|^2\bigr) =A-4\uppi
f_{X,j}B.
\end{eqnarray*}

Next, use the fact that $ | |z_1|^2-|z_2|^2 |\le
|z_1-z_2|^2+2|z_1-z_2| |z_2|$,
for any complex numbers $z_1, z_2$,
and that $|f_{X\zeta,j}|=|A_{X,j}|/(2\uppi)^2\le f_{X,j}^{1/2}$, to obtain
\begin{eqnarray*}
|A-4\uppi f_{X,j}B| 
&\le& |A|+4\uppi
f_{X,j}|B|
\\
&\le& \bigl(E|w_{X,j}|^2-f_{X,j}
\bigr)^2+2f_{X,j} \bigl|E|w_{X,j}|^2-f_{X,j}
\bigr|
\\
&&{} +4\uppi f_{X,j} \bigl|E[w_{X,j}\overline{w_{\zeta
,j}}]-f_{X\zeta,j}
\bigr|^2 +8\uppi f_{X,j}^{3/2} \bigl|E[w_{X,j}
\overline{w_{\zeta,j}}]-f_{X\zeta
,j} \bigr|,
\end{eqnarray*}
which shows that $s_{n,1}^\prime\le C s_{n,1}$ and
completes the proof of the corollary.
\end{pf}

\begin{pf*}{Proof of Lemma \ref{lemlap}} The proof uses Proposition
\ref{propf-mult}.
We shall prove (\ref{es-con0>}) and (\ref{es-con0?}).
These two facts together imply (\ref{esacon}) in a routine fashion. %

\textit{Proof of} (\ref{es-con0>}). By (\ref{esxsz}), $R_n$ is
like the r.v. in the l.h.s. of (\ref{eri-co2}) with
$c_{n,j}= b_{n,j}/f_{X,j}$. Thus, $\operatorname {Var}(R_n)\le s_{n,1}+s_{n,2}$,
where $s_{n,k}, k=1,2$
are the same in (\ref{eri-co2}) with
$c_{n,j}\equiv b_{n,j}/f_{X,j}$.
It thus suffices to show that the sum $s_{n,1}+s_{n,2}$ is bounded from
the above by
the upper bounds given in
(\ref{es-con0>}). 

Recall Proposition \ref{propf-mult}(iii).
The spectral density $f_X$ satisfies (\ref{effNew}), whereas the
cross-spectral density $f_{X\zeta}(u)=(2\uppi)^{-1}A_X(u)$ has the
property $|f_{X\zeta}(u)|\le C|u|^{-d}$, $|\dot{f}_{X\zeta}(u)|\le
C|u|^{-1-d}$, $u\in\Pi$. Therefore, they satisfy conditions of this
proposition, and hence
%
\begin{eqnarray}
\label{efisrta1} \bigl|E|w_{X,j}|^2-f_{X,j} \bigr| &\le&
C|u_j|^{-2d}j^{-1}\log(1+j),
\nonumber
\\[-8pt]
\\[-8pt]
\nonumber
\bigl|E[w_{X,j}\overline{w_{\zeta,j}}]-f_{X\zeta,j} \bigr| &\le&
C|u_j|^{-d}j^{-1}\log (1+j),
\end{eqnarray}
where $C$ does not depend on $j$ and
$n$. Since, by (\ref{effNew}), $1/f_{X,j}(u)\le Cu_j^{2d}$, these
bounds yield
$ 
s_{n,1}\le C \sum_{j=1}^{\nu} b_{n,j}^2
(j^{-1}\log j).
$ 
This bound, the Cauchy--Schwarz inequality, and the fact $\sum_{j\ge1}
j^{-2}\log^2 j<\iny,$
imply
\[
s_{n,1}\le C b_n^2\log(n)\sum
_{j=1}^{\nu} j^{-1}\le C b^2_n
\log^2(n)\quad \mbox{and}\quad 
s_{n,1}\le C
b_n\sum_{j=1}^{\nu}
|b_{n,j}| \bigl(j^{-1}\log j\bigr) \le
Cb_nB_n.
\]
This proves that $s_{n,1}$ satisfies both bounds of (\ref{es-con0>}).

Next, again by Proposition \ref{propf-mult}(iii), for all $1\le
k<j\le
\nu$,
\begin{eqnarray*}
\bigl|E[w_{X,j}\overline{w_{X,k}}] \bigr|&\le& C\bigl(u_j^{-2d}+u_k^{-2d}
\bigr)j^{-1}\log j,\qquad \bigl|E[w_{X, j}\overline{w_{\zeta,k}}] \bigr|
\le C\bigl(u_j^{-d}+u_k^{-d}
\bigr)j^{-1}\log j.
\end{eqnarray*}
By (\ref{effNew}),
\begin{eqnarray*}
(f_jf_k)^{-1}\bigl(u_j^{-2d}+u_k^{-2d}
\bigr)^2&\le& C(u_ju_k)^{2d}
\bigl(u_j^{-4d}+u_k^{-4d}\bigr)\le
C(j/k)^{2|d|},
\\
f_j^{-1}\bigl(u_j^{-d}+u_k^{-d}
\bigr)^2&\le& Cu_j^{2d}\bigl(u_j^{-2d}+u_k^{-2d}
\bigr)\le C(j/k)^{2|d|}.
\end{eqnarray*}
These facts together imply
%
\begin{equation}
\label{edaw1} s_{n,2} \le C\sum_{1\le k<j\le\nu}
|b_{n,j}b_{n,k}| \biggl(\frac{j}{k}
\biggr)^{2|d|}\frac{\log^2 j}{j^2} . 
\end{equation}

Bound $|b_{n,j}b_{n,k}|$ by $b_n^2$ to obtain
$ s_{n,2}\le C b_n^2\log^2(n)\sum_{1\le k<j\le\nu} k^{-2|d|}j^{2|d|-2}
\le\break C b_n^2 \log^3(n),$
which implies the first estimate of (\ref{es-con0>}).
Next, bound $|b_{n,j}|$ by $b_n $ in (\ref{edaw1}), to obtain
\begin{eqnarray*}
s_{n,2}&\le& Cb_n \sum_{1\le k<j\le\nu}
|b_{n,k}| \frac{\log^2 j}{k^{2|d|}j^{2-2|d|}} \le C b_n \sum
_{1\le k\le\nu} |b_{n,k}| \frac{\log^2 k}{k}
\\
&\le& Cb_n \biggl(\sum_{1\le k\le\nu}
b_{n,k}^2 \biggr)^{1/2} \biggl(\sum
_{1\le k\le\nu}\frac{\log^4 k}{k^2} \biggr)^{1/2} \le
Cb_n B_n,
\end{eqnarray*}
which establishes the second bound of
(\ref{es-con0>}).

To show (\ref{es-con0?}), recall
that $f_{X,j}E|w_{\zeta, j}|^2/f_{\zeta,j}=f_{X,j}$.
Therefore,
\[
ER_n=\sum_{j=1}^{\nu}
\frac{b_{n,j}}{f_{X,j}} \biggl(E|w_{X,j}|^2- \frac{f_{X,j}}{f_{\zeta,j}}E|w_{\zeta,
j}|^2
\biggr) = \sum_{j=1}^{\nu}\frac{b_{n,j}}{f_{X,j}}
\bigl(E|w_{X,j}|^2- f_{X,j} \bigr).
\]
Then, by (\ref{efisrta1}) and (\ref{effNew}),
\begin{eqnarray*}
|ER_n|\le C\sum_{j=1}^{\nu}
\frac
{|b_{n,j}|}{f_{X,j}}u_j^{-2d}j^{-1}\log j \le C\sum
_{j=1}^{\nu}|b_{n,j}|j^{-1}
\log j 
\le Cb_n\log^2(n),
\end{eqnarray*}
which implies the first bound in (\ref{es-con0?}).\eject

To establish the second bound, let $K=(B_n/b_n)^{1/2}$.
Because of (\ref{es-con1}), $K\to\infty$, $b_nK=(b_n/B_n)^{1/2}B_n=\mathrm{o}(B_n)$.
Thus,
%
\begin{eqnarray}
\label{earg} |ER_n|&\le&C \Biggl(\sum_{j=1}^{K-1}|b_{n,j}|j^{-1}
\log j+\sum_{j=K}^{\nu} |b_{n,j}|j^{-1}
\log j \Biggr)
\nonumber
\\[-8pt]
\\[-8pt]
\nonumber
&\le& C \Biggl\{b_nK+ \Biggl(\sum_{j=K}^{\nu}b_{n,j}^2
\Biggr)^{1/2} \Biggl(\sum_{j=K}^\infty
j^{-2}\log^2 j \Biggr)^{1/2} \Biggr\}
=\mathrm{o}(B_n).
\end{eqnarray}
This completes proof of the second estimate
in (\ref{es-con0?}).
\end{pf*}

Now we establish the asymptotic normality of the weighted quadratic
forms $S_{n,\zeta}$.
The CLT for quadratic forms
in i.i.d. r.v.'s is well investigated; see
Guttorp and Lockhart \cite{GutLoc88}. The following theorem summarizes a useful criterion
for asymptotic normality, given in Theorem 2.1 in
Bhansali \textit{et al.} \cite{BhaGirKok07N1}. Let $C_n=\{c_{n,ts}, t,s=1,\ldots
,n\}$ be
a symmetric $n\times n$ matrix of real numbers $c_{n,ts}$, and define
the quadratic form
\[
\cQ_{n}:=\sum_{t,s=1}^n
c_{n,ts}\zeta_t\zeta_s.
\]
Let $\|C_n\|:=(\sum_{t,s=1}^n c_{n,ts}^2)^{1/2}$ and
$\|C_n\|_{\mathrm{sp}}:=\max_{\|x\|=1}\|C_n x\|$
denote Euclidean and spectral norms, respectively, of $C_n$.
%
\begin{thmm}
\label{tBGK} Suppose $\zeta_j\sim \operatorname{IID}(0,1)$ and $E\zeta_0^4<\infty$.
Then
%
\begin{equation}
\label{ebb+} \frac{\|C_n\|_{\mathrm{sp}}}{\|C_n\|}\rightarrow0
\end{equation}
implies
$ 
(\operatorname{ Var}(\cQ_n))^{-1/2}(\cQ_n-E\cQ_n)\to_D \cN(0,1).
$ 
In addition, if 
$\sum_{t=1}^n c_{n;tt}^2=\mathrm{o}(\|C_n\|^2)$, then $\operatorname{ Var}(\cQ_n)\sim
2\|C_n\|^2$. Furthermore, in this case, if $E\zeta_0^4<\infty$ is
replaced by $E|\zeta_0|^{2+\delta}<\infty$, for some $ \delta>0$, then
$(2\|C_n\|^2)^{-1/2}(\cQ_n-E\cQ_n)\to_D \cN(0,1).$
\end{thmm}
%

Next lemma derives asymptotic distribution of the sum
$S_{n,\zeta}$ of (\ref{esnz}). Its proof uses Theorem~\ref{tBGK}
and some ideas of the proof of Theorem 2, Robinson \cite{ShaWu07N2}.
%
\begin{lem}\label{lemlasymp} Suppose $\zeta_j\sim \operatorname{IID}(0,1)$,
$E\zeta_0^{4}<\infty$, and $b_{n,j}$
satisfy (\ref{es-con1}).
Then
%
\begin{eqnarray}
\label{es-mean} &&\mathrm{(a)}\quad ES_{n,\zeta}= \sum_{j=1}^{\nu}
b_{n,j},\qquad 
\mathrm{(b)}\quad \operatorname{Var} (S_{n,\zeta})=
q_n^2,
\nonumber
\\[-8pt]
\\[-8pt]
\nonumber
&&\mathrm{(c)}\quad q_n^{-1}(S_{n,\zeta}-ES_{n,\zeta})
\to_D \cN(0,1).
\end{eqnarray}
Moreover,
%
\begin{equation}
\label{ecum4} q_n^2\ge\min \bigl(1,\operatorname{ Var}
\bigl(\zeta_0^2\bigr)/2 \bigr)B_n^2.
\end{equation}
\end{lem}
\begin{pf} Let $c_n(t):= n^{-1}\sum_{j=1}^{\nu}b_{n,j} \cos(tu_j),
t=1,2,\ldots.$ Note that
\[
S_{n, \zeta}=\frac{1}{n}\sum_{t,s=1}^n
\sum_{j=1}^{\nu} \mathrm{e}^{\i(t-s)u_j}
b_{n,j}\zeta_s \zeta_t =\sum
_{t,s=1}^n c_n(t-s)\zeta_s
\zeta_t.
\]
The matrix
$C_n=(c_n(t-s))_{t,s=1,\ldots,n}$ is a symmetric $n\times n$
matrix with real entries. Hence, (\ref{es-mean})(a) and
(\ref{es-mean})(b) follow because $\zeta_j$'s are $\operatorname{IID}(0,1)$. For the
same reason, and because $\operatorname{ Var}(\zeta_0^2)-2= E\zeta_0^4-3=\operatorname{
Cum}_4(\zeta_0)$, and $ c_n(0)= n^{-1}\sum_{j=1}^{\nu}b_{n,j},$
%
\begin{eqnarray}
\label{e2fa} \operatorname{ Var} (S_{n,\zeta})&=&2\sum
_{s,t=1: t\ne s}^n c^2_n(t-s) +
\operatorname{ Var}\bigl(\zeta_0^2\bigr)\sum
_{t=1}^n c^2_n(t-t)
\nonumber
\\[-8pt]
\\[-8pt]
\nonumber
&=&2\|C_n\|^2+\operatorname{ Cum}_4(
\zeta_0)n^{-1}\Biggl(\sum_{j=1}^{\nu
}b_{n,j}
\Biggr)^2\ge 
\min\bigl(2, \operatorname{ Var}\bigl(
\zeta_0^2\bigr)\bigr)\|C_n
\|^2.
\end{eqnarray}
We shall show below that
%
\begin{eqnarray}
\label{eWWW+} \mathrm{(a)}\quad \|C_n\|^2 =2^{-1}B_n^2,\qquad
\mathrm{(b)}\quad  \|C_n\|_{\mathrm{sp}}=\mathrm{o}(\|C_n
\|).
\end{eqnarray}
Then (\ref{e2fa}) and (\ref{eWWW+})(a) imply (\ref{ecum4}), 
whereas by Theorem \ref{tBGK}, (\ref{eWWW+})(b)
implies
\begin{eqnarray*}
\bigl(\operatorname {Var}(S_{n,\zeta})\bigr)^{-1/2}\bigl(S_{n,\zeta}-E[S_{n,\zeta}]
\bigr) &\to_D& \cN(0,1),\\
 \operatorname {Var} (S_{n,\zeta})&=&B_n^2+
\operatorname{ Cum}_4(\zeta_0)n^{-1}\Biggl(\sum
_{j=1}^{\nu}b_{n,j}
\Biggr)^2,
\end{eqnarray*}
which proves (\ref{es-mean})(c). 
It remains to show (\ref{eWWW+}).


To prove (\ref{eWWW+})(a), recall that for all $1\le j,k\le m,
j+k<n$ and $ a, b\in\R$,
%
\begin{equation}
\sum_{t=1}^n\cos(tu_j+a)
\cos(tu_k+b)=\frac{n}{2}\cos(a-b)I(j=k) . \label{e00}
\end{equation}

This fact and the definition of $c_n(t)$ imply (\ref{eWWW+})(a),
because
\[
\label{e009} \|C_n\|^2 =\sum
_{t,s=1}^n c_n^2(t-s) =
n^{-2}\sum_{j,
k=1, j+k<n}^{\nu}b_{n,j}b_{n,k}
\sum_{s,t=1}^n \cos\bigl((t-s)u_j
\bigr)\cos \bigl((t-s)u_k\bigr)
= 2^{-1} B_n^2.
\]

To establish (\ref{eWWW+})(b), note that by (\ref{e00}),
$ 
\sum_{t=1}^n c_n(t-s)c_n(t-v)
=
(2n)^{-1}\times \sum_{j=1}^{\nu} b_{n,j}^2\cos((s-v)u_j).
$ 
Hence, for any $x\in\R^n$, such that $\|x\|=1$,
\begin{eqnarray*}
\label{amam} \|C_n x\|^2 &=& \sum
_{t=1}^n \Biggl(\sum_{s=1}^n
c_n(t-s)x_s \Biggr)^2 = \sum
_{s, v=1}^n x_sx_v \Biggl(
\sum_{t=1}^nc_n(t-s)c_n(t-v)
\Biggr)
\\
&=& \frac{1}{2n}\sum_{j=1}^{\nu}
b_{n,j}^2\sum_{s,v=1}^n
\cos \bigl((s-v)u_j\bigr)x_sx_v 
\le\frac{1}{2n} b_n^2\sum
_{j=1}^n \Biggl|\sum_{s=1}^ne^{\i su_j}x_s
\Biggr|^2. 
\end{eqnarray*}
Expand the last quadratic and use the fact $\sum_{j=1}^ne^{\i
(t-s)u_j}=nI(t=s)$, to obtain
\begin{eqnarray*}
\|C_n x\|^2 &\le& \frac{1}{2} b_n^2
\sum_{t=s=1}^nx_t^2=
\frac{1}{2} b_n^2\|x\|^2,\qquad
\|C_n\|_{\mathrm{sp}}\le(1/\sqrt2){b_n}.
\end{eqnarray*}
Since $b_n=\mathrm{o}(B_n)$, and
$B_n=\sqrt2 \|C_n\|$ by (\ref{eWWW+})(a), this proves
(\ref{eWWW+})(b), and also completes the proof of the lemma.
\end{pf}

\section{A general case of sums of weighted periodogram}\label {sgeneral}
We now focus on the sums $Q_{n,X}$ of (\ref{esum11}).
Bartlett approximation 
$I_{X,j}\sim f_{X,j} (I_{\zeta,j}/f_{\zeta,j} )$
suggests to approximate $Q_{n,X}$ by the sum
\[
Q_{n,\zeta}:= \sum_{j=1}^{ \nu}
(b_{n,j} f_{X,j}) \biggl(\frac{I_{\zeta
,j}}{f_{\zeta,j}}\biggr)= \sum
_{j=1}^{ \nu} b_{n,j}
f_{X,j}(2\uppi)I_{\zeta,j}.
\]
%

In Theorem \ref{ts-main} above, 
$f_X$
can be unbounded at $0$, but differentiable on $(0,\uppi)$.
Then the asymptotic normality
of the sums
$S_{n,X}=\sum_{j=1}^{ \nu} b_{n,j} (I_{X,j}/f_{X,j} )$
holds
under 
(\ref{es-con1}).

Now we turn to the case when $f_X$ is continuous on $\Pi$ and satisfies
%
\begin{eqnarray}
\label{eff-bounded} 0<C_1\le f_X(u)\le
C_2<\infty, \qquad u\in\Pi\ (\exists 0<C_1,C_2<
\infty).
\end{eqnarray}
%

Theorems \ref{tcontinuous} and \ref{thmbias-t}
below show that under 
(\ref{es-con1}),
continuity of $f_X$, or more precisely,
continuity of the transfer function $A_X$, suffices for
asymptotic normality of the centered sums
$Q_{n,X}-EQ_{n,X}$ and for obtaining an upper bound on the variance
$\operatorname {Var}(Q_{n,X})$,
whereas satisfactory asymptotics of $EQ_{n,X}$ requires $f_X$ to be
Lipshitz($\b$), $\b>1/2$.

By Lemma \ref{lemlasymp}, 
$EQ_{n,\zeta}=\sum_{j=1}^{\nu}b_{n,j}f_{X,j}$ and $\operatorname{
Var}(Q_{n,\zeta})=v_n^2,$ where
\[
\label{eqbf} v_n^2:=\sum_{j=1}^{\nu}(b_{n,j}f_{X,j})^2+
\operatorname{ Cum}_4(\zeta_0)\frac{1}{n} \Biggl(
\sum_{j=1}^{\nu
}b_{n,j}f_{X,j}
\Biggr)^2.
\]
%
Let $ 
b_{f,n}=\max_{j=1,\ldots,\nu}|b_{n,j}|f_{X,j}$ and $ B_{f,n}^2=
\sum_{j=1}^{\nu}(b_{n,j}f_{X,j})^2.
$ Similarly as in (\ref{ebab14}), one can show that
for some $C_1,C_2>0$,
%
\begin{eqnarray}
\label{eqbf-1} C_1B^2_{f,n}\le
v_n^2 \le C_2B_{f,n}^2\quad
\mbox{and}\quad C_1B^2_{n}\le
v_n^2\le C_2B_{n}^2,
\qquad\mbox{under } (\ref{eff-bounded}).
\end{eqnarray}

The folowing theorem describes the asymptotic behavior of bias,
variance, and asymptotic normality of $Q_{n,X}$ when $f_X$ is
continuous and bounded.
%
\begin{thmm}\label{tcontinuous} Suppose the linear process $\{X_j,
j\in\Z\}$
of (\ref{elinear}) is such that $E\zeta^4_0$ $<\infty$,
and the real weights $b_{n,j}$'s satisfy (\ref{es-con1}).

In addition, if $f_X$ satisfies (\ref{eff-bounded}) and $A_X\in\cC
(\Pi)$,
then
%
\begin{eqnarray}
\label{econtin-var} \operatorname{ Var}(Q_{n,X})= v_n^2+o
\bigl(v_n^2\bigr), \qquad v_n^{-1}(Q_{n,X}-EQ_{n,X})
\to_D \cN(0,1).
\end{eqnarray}
In addition, if $f_X\in\Lambda_\beta(\Pi)$, with
$\beta>1/2$, then
%
\begin{equation}
\label{ec-mean} E Q_{n,X} = \sum_{j=1}^{\nu}
b_{n,j}f_{X,j}+\mathrm{o}(v_n), \qquad v_n^{-1}
\Biggl(Q_{n,X}-\sum_{j=1}^{\nu}
b_{n,j}f_{X,j}\Biggr)\to_D \cN(0,1).
\end{equation}
\end{thmm}

The next theorem covers the 
case when the $f_X$
is not bounded in the neighborhood of $0$, that is, $d>0$, or is not
bounded away from $0$,
that is, $d<0$.
Then the second bound of (\ref{eqbf-1}) does not hold.
Assumption (\ref{es-con1}) now has to be formulated
using the weights
$b_{n,j}f_{X,j}$ and 
we need to impose some additional smoothness
conditions on $A_X$ in a small neighborhood of~$0$. We assume
that
$A_X$ can be factored into a product $A_X=hG$ of a differentiable function
$h$, which may have a pole at~$0$, and a continuous bounded function $G$.
In particular, if $A_X$ satisfies (\ref{eaaa}), we take
$G\equiv1$.
%
\begin{thmm}\label{tcontinuous1} Suppose $\{X_j, j\in\Z\}$ is the
linear process
(\ref{elinear}) with $E\zeta^4_0$ $<\infty$.
Assume that $f_X$ satisfies (\ref{effNew}) with 
$|d|<1/2$, the transfer function $A_X$
can be factored as $A_X=hG$, where $G$
is continuous
and bounded away from $0$
and $\infty$,
and $h$ is differentiable having derivative $\d h$ and satisfying
%
\begin{equation}
\label{ehhmark} C_1|u|^{-d}\le\bigl|h(u)\bigr|\le
C_2|u|^{-d},\qquad  \bigl|\dot{h}(u)\bigr|\le C|u|^{-1-d},\qquad  0<|u|\le
\uppi,
\end{equation}
for some $0<C, C_1, C_2<\iny$. Then, for any real weights $b_{n,j}$'s
satisfying
%
\begin{equation}
\label{eLF-1} \frac{b_{f,n}}{B_{f,n}} \equiv\frac{\max_{j=1,\ldots, \nu
}|b_{n,j}f_{X,j}|}{(\sum_{j=1}^{\nu}(b_{n,j}f_{X,j})^2)^{1/2}} \to0,
\end{equation}
%
(\ref{econtin-var}) continues to hold.

If, in addition, $G\in\Lambda_\beta(\Pi)$, with $\beta>1/2$,
then also (\ref{ec-mean}) holds.
\end{thmm}
\begin{pf*}{Proofs of Theorems~\ref{tcontinuous} and~\ref{tcontinuous1}}
The proof of both theorems follows from Lemmas~\ref{lemlasymp} and~\ref{lemlLip}. The latter lemma will be proved
shortly.

Let $ r_n:=Q_{n,X}-Q_{n,\zeta}-E[Q_{n,X}-EQ_{n,\zeta}]. $ In Lemmas~\ref{lemlLip}(i)
and~\ref{lemlLip}(ii), 
it is shown that $E
r_n^2=\mathrm{o}(v_n^2)$ under the assumptions of Theorems~\ref{tcontinuous}
and~\ref{tcontinuous1}. Therefore, the claim (\ref{econtin-var})
made in these two theorems 
follows, noticing that, by Lemma \ref{lemlasymp}, under assumption
(\ref{eLF-1}), $v_n^{-1}(Q_{n,X}-Q_{n,\zeta})\to_D\mathcal{
N}(0,1)$. The
second claim (\ref{ec-mean}) of these theorems is shown in
(\ref{emeanmean}) of Theorem \ref{thmbias-t} below.
\end{pf*}

Lemma \ref{lemlLip} below shows that the order of
approximation of $Q_{n,X}-EQ_{n,X}$ by $Q_{n,\zeta}-EQ_{n,\zeta}$ is
determined by the smoothness of the transfer function $A_X$.
For example, by Lemma \ref{lemlLip}(i), if $A_X$ is a bounded
continuous function, then
%
\begin{equation}
\label{eaqpa} Q_{n,X}-EQ_{n,X}= Q_{n,\zeta}-EQ_{n,\zeta}+\mathrm{o}_p(v_n).
\end{equation}
%
If, in addition, $A_X$
has a bounded derivative, then the order
improves to $\mathrm{o}_p(n^{-1/2}\log(n) v_n)$ without requiring any
additional
assumptions on $b_{n,j}$. Lemma \ref{lemlLip}(ii) shows that if $A_X$
is discontinuous at $0$,
then approximation (\ref{eaqpa}) is valid under additional
regularity behavior of $A_X$ in a neighborhood of $0$, as long as
the weights $b_{n,j}$ satisfy (\ref{eLF-1}).

To state the lemma, we need the following notation. For a
complex valued function $h(u), u\in\Pi$, define
\begin{eqnarray*}
\label{edelta} \varepsilon_{n,h}&:=& n^{-1}
\log^2(n),\qquad \mbox{$h\in\Lambda_1[\Pi]$,}
\\
&:=&n^{-\beta},\qquad \mbox{$h\in\Lambda_\beta[\Pi]$, $0<\beta <1$},
\\
&:=&\delta_n, \delta_n\to0, \qquad\mbox{$h\in\mathcal{ C}[
\Pi]$}.
\end{eqnarray*}

\begin{lem}\label{lemlLip} Assume that $\{X_j\}$ is as in (\ref
{elinear}) and
$E\zeta^4_0<\infty$.
Then 
the following hold.
\begin{longlist}[(ii)]
\item[(i)]
If $A_X\in\Lambda_\beta[\Pi]$, $0<\beta\le1$, or $A_X\in\mathcal
{ C}[\Pi]$, then
%
\begin{equation}
\label{esacLip??} Er^2_n \le C\varepsilon_{n,A_X}
B_n^2=\mathrm{o}\bigl(v_n^2\bigr).
\end{equation}

\item[(ii)]
If $A_X=hG$, where $h$ satisfies (\ref{ehhmark})
and either $G\in\mathcal{ C}(\Pi)$ or $G\in\Lambda_\beta(\Pi)$,
$0< \beta\le1$, then
%
\begin{eqnarray}
\label{esac-A} Er^2_n&\le& C \bigl(\min
\bigl(b_{f,n}^2 \log^3 n , b_{f,n}
B_{f,n}\bigr)+ \varepsilon_{n,G} B_{f,n}^2
\bigr),
\nonumber
\\[-8pt]
\\[-8pt]
\nonumber
&\le& C\min\bigl(b_{f,n}^2 \log^3 n ,
b_{f,n} B_{f,n}\bigr), \qquad G\in \Lambda_1(\Pi).
\end{eqnarray}
If, in addition, (\ref{eLF-1}) holds,
then
%
\begin{equation}
\label{erabo} Er^2_n=\mathrm{o}\bigl(v_n^2
\bigr).
\end{equation}
\end{longlist}
\end{lem}

\begin{pf} Rewrite $r_n=D_n -ED_n$, where\vspace*{1pt}
$D_n=Q_{n,X}-Q_{n,\zeta}=\sum_{j=1}^{\nu}b_{n,j}\{
I_{X,j}-(f_{X,j}/f_{\zeta,j})\times I_{\zeta,j}\}.
$
Let $t_{n,i}, i=1,2$ denote the $s_{n,i}, i=1,2,$ of Corollary
\ref{ccor-b} with
$c_{n,j}\equiv b_{n,j}$. By Corollary~\ref{ccor-b},
%
\begin{equation}
\label{evart}
\operatorname{ Var}(D_n) \le C(t_{n,1}+t_{n,2}).
\end{equation}
%

\textit{Proof of} (i). Arguing as in the proof of Lemma \ref{lemlap},
one can show that
%
\begin{equation}
\label{es-next} E\bigl(r_n^2\bigr)\le C
\varepsilon_{n,A_X} B_n^2,
\end{equation}
which, in view of (\ref{eqbf-1}), proves (\ref{esacLip??}). We need
to verify (\ref{es-next}) 
in the following three cases.

Case (1). $A_X\in\Lambda_1[\Pi]$. Then, by Proposition
\ref{propf-mult}(i),
\begin{eqnarray*}
\bigl|E|w_{X,j}|^2-f_{X,j}\bigr|
\vee\bigl|E[w_{X,j}\overline{w_{\zeta
,j}}]-f_{X\zeta,j}\bigr|&\le&
Cn^{-1}\log n,
\\
\bigl|E[w_{X,j}\overline{w_{X,k}}]\bigr|\vee \bigl|E[w_{X,j}
\overline{w_{\zeta,k}}]\bigr|&\le& Cn^{-1}\log n, \qquad 1\le k<j\le\nu.
\end{eqnarray*}
Therefore, $t_{n,1} \le Cn^{-1}\log
n\sum_{j=1}^{ \nu}b_{n,j}^2= Cn^{-1}\log n B^2_n,$ and
\[
t_{n,2}\le Cn^{-2}\log^2 n\sum
_{1\le k<j\le\nu}|b_{n,j}b_{n,k}|
\le Cn^{-1}\log^2 nB^2_n,
\]
which proves (\ref{es-next}).

Case (2). $A_X\in\Lambda_\beta[\Pi]$, $0<\beta<1$. Then by
Proposition \ref{propf-mult}(i),
\begin{eqnarray*}
\bigl|E|w_{X,j}|^2-f_{X,j}\bigr|
\vee\bigl|E[w_{X,j}\overline{w_{\zeta
,j}}]-f_{X\zeta,j}\bigr|&\le&
Cn^{-\beta},
\\
\bigl|E[w_{X,j}\overline{w_{X,k}}]\bigr|\vee \bigl|E[w_{X,j}
\overline{w_{\zeta,k}}]\bigr|&\le& Cn^{-\beta}\ell_n(\beta;j-k),\qquad
k<j.
\end{eqnarray*}
Note that for $ 1\le k<j\le\nu<n/2$, $j-k\le n-j+k$, and
hence bound
\[
\ell_n(\beta;j-k) \le C \frac{\log(2+ j-k)}{(2+j-k)^{1-\beta}},\qquad
\bigl(n^{-\beta}\ell_n(\beta;j-k)
\bigr)^2\le C\frac{\log^2
(2+j-k)}{n^{\beta}(2+j- k)^{2-\beta}}. 
\]
Apply this fact, to obtain, that for $0<\beta<1$,
$t_{n,1}\le Cn^{-\beta} \sum_{j=1}^{ \nu}b_{n,j}^2= Cn^{-\beta}
B^2_n$,
\begin{eqnarray*}
t_{n,2}&\le&C \sum
_{1\le k<j\le\nu}|b_{n,j}b_{n,k}|
\bigl(n^{-\beta}\ell_n(\beta;j-k)\bigr)^2
\\
&\le& Cn^{-\beta} \sum_{1\le k<j\le\nu}|b_{n,j}b_{n,k}|
\frac{\log^2 (2+j-k)}{(2+j-k)^{2-\beta}}
\le Cn^{-\beta}B_n^2,
\end{eqnarray*}
which proves (\ref{es-next}).

Case (3). $A_X\in\mathcal{ C}[\Pi]$. By Proposition
\ref{propf-mult}(ii), for any $0<\varepsilon<1/2$,
\begin{eqnarray*}
\bigl|E|w_{X,j}|^2-f_{X,j}\bigr|
\vee\bigl|E[w_{X,j}\overline{w_{\zeta
,j}}]-f_{X\zeta,j} \bigr| &\le& C
\delta_n,
\\
\bigl|E[w_{X,j}\overline{w_{X,k}}]\bigr|\vee\bigl|E[w_{X,j}
\overline {w_{\zeta,k}}] \bigr| &\le & C\delta_n 
\ell_n(\varepsilon;j-k),\qquad  k<j,
\end{eqnarray*}
with some $\delta_n\to0$, that does not depend on $k, j$
and $n$,
and (\ref{es-next}) follows by the same argument as in the case
(2) above. This completes the proof of (i) of the lemma.

\textit{Proof of} (ii). First, we prove
(\ref{esac-A}). As above, for that we need to bound
$t_{n,1}$ and $t_{n,2}$ of (\ref{evart}).
Recall that $f_{X}=|A_X|^2/(2\uppi),$ $f_{X\zeta}=A_X/(2\uppi)$,
$A_X=h(u)G(u)$, where $h$ satisfies (\ref{ehhmark}), which
together with (\ref{effNew}) implies that $G$ is bounded away
from infinity and zero. For $1\le k\le j\le\nu,$ define
\begin{eqnarray*}
\tilde r_{n,jk}&:=&0, \qquad G\in\Lambda_1(\Pi),
\\
&:=&n^{-\beta} \frac{\log(2+j-k)}{(2+j-k)^{1-\beta}},\qquad G\in\Lambda_\beta(\Pi), 0<
\beta<1,
\\
&:=&\delta_n \frac{\log(2+j-k)}{(2+j-k)^{1-\varepsilon}},\qquad  G\in\mathcal{ C}(\Pi), 0<
\varepsilon<1/2, \del_n\to0.
\end{eqnarray*}
%

By Proposition \ref{propf-mult}(iv), for $1\le k\le j$,
\begin{eqnarray*}
 \bigl|E[w_{X,j}\overline{w_{X,k}}]-f_{X,j}I(j=k)\bigr|
&\le & C\bigl\{\bigl(u_k^{-2d}+u_j^{-2d}
\bigr)j^{-1} \log j + \bigl(u_k^{-2d}\wedge
u_j^{-2d} \bigr) \tilde r_{n,jk}\bigr\}
\\
\bigl|E[w_{X,j}\overline{w_{\zeta,k}}]-f_{X\zeta,j}I(j=k)\bigr|
&\le& C\bigl\{\bigl(u_k^{-d}+u_j^{-d}
\bigr)j^{-1} \log j + \bigl(u_k^{-d}\wedge
u_j^{-d}\bigr) \tilde r_{n,jk}\bigr\}.
\end{eqnarray*}
Since
$f_X=|A_X|^2/(2\uppi)=|hG|^2/(2\uppi)$, assumptions on $h$ and $G$ here
imply that for all $u\in\Pi$,
\begin{eqnarray*}
f_X(u)\le C|u|^{-2d},\qquad f_X^{-1}(u)
\le C|u|^{2d},\qquad \bigl|f_{X\zeta}(u)\bigr|\le C|u|^{-d},\qquad
\bigl|f^{-1}_{X\zeta}(u)\bigr|\le C|u|^{d}.
\end{eqnarray*}
Therefore, for $ 1\le k\le j$,
\begin{eqnarray*}
(f_{X,j}f_{X,k})^{-1}\bigl(u_k^{-2d}+u_j^{-2d}
\bigr)^2&\le &C|j/k|^{2|d|}, \qquad (f_{X,j}f_{X,k})^{-1}
\bigl(u_k^{-2d}\wedge u_j^{-2d}
\bigr)^2\le C,
\\
(f_{X,j})^{-1}\bigl(u_k^{-d}+u_j^{-d}
\bigr)^2&\le& C|j/k|^{2|d|}, \qquad (f_{X,j})^{-1}
\bigl(u_k^{-d}\wedge u_j^{-d}
\bigr)^2\le C.
\end{eqnarray*}
%

Recall the bound (\ref{evart}).
It suffices to show that $t_{n,1}+t_{n,2}$ can be bounded above by
the r.h.s. of (\ref{esac-A}). The above bounds readily yield that
\begin{eqnarray*}
t_{n,1}&\le& C\sum_{j=1}^{\nu}
(b_{n,j}f_{X,j})^2 \bigl(j^{-1}\log j +
\tilde r_{n,jj}\bigr),
\\
t_{n,2}&\le& C\sum_{1\le k<j\le\nu}
|b_{n,j}f_{X,j}| |b_{n,k}f_{X,k}| \biggl(
\biggl(\frac{j}{k}\biggr)^{2|d|}\frac{\log^2
j}{j^2}+\tilde
r^2_{n,jk} \biggr).
\end{eqnarray*}
The arguments analogous to one used in evaluating $s_{n,1}$ and $s_{n,2}$
in Lemma \ref{lemlap} yield
\begin{eqnarray*}
&&\sum_{j=1}^{\nu} (b_{n,j}f_{X,j})^2
\frac{\log
j}{j}+\sum_{1\le k< j\le\nu} |b_{n,j}f_{X,j}|
|b_{n,k}f_{X,k}|\biggl(\frac{j}{k}
\biggr)^{2|d|}\frac{\log^2 j}{j^2 }
\\
&&\quad \le C\min \bigl( b^2_{f,n}\log^3(n),
b_{f,n} B_{f,n} \bigr),
\\
&& \sum_{j=1}^{\nu}
(b_{n,j}f_{X,j})^2\tilde r_{n,jk}+ \sum
_{1\le k< j\le\nu} |b_{n,j}f_{X,j}|
|b_{n,k}f_{X,k}| \tilde r_{n,jk}^2\le C
\varepsilon_{n,G} B_{f,n}^2.
\end{eqnarray*}
Therefore,
$ 
t_{n,1}+t_{n,2}\le C (\min (
b^2_{f,n}\log^3(n), b_{f,n} B_{f,n} )+ \varepsilon_{n,G}
B_{f,n}^2 ),$
which proves (\ref{esac-A}).

Observe that $\varepsilon_{n,G}\to0$. Therefore, (\ref{esac-A}),
(\ref{eLF-1}) and (\ref{eqbf-1}) imply (\ref{erabo}). This
completes the proof of the lemma.
\end{pf}

As seen above, proving CLT for
$v_n^{-1}(Q_{n,X}-\sum_{j=1}^{\nu}b_{n,j}f_{X,j})$ requires some
smoothness of the spectral density $f_X$ and the transfer function
$A_X$. Conditions on $A_X$ can be relaxed if one wishes to
establish only an upper bound for the mean square error of the
estimator $Q_{n,X}$ of $\sum_{j=1}^{\nu}b_{n,j}f_{X,j}$ as is shown
in the next theorem. The results of Theorem \ref{thmbias-t} also
remain valid for $\nu=[n/2]$.
%
\begin{thmm}\label{thmbias-t}
Let $\{X_j\}$ be as in (\ref{elinear}) with
$E\zeta^4_0<\infty$ and $f_X$ satisfying (\ref{effNew}).
\begin{enumerate}[(ii)]
\item[(i)] Then
%
\begin{equation}
\label{esumsum} E(Q_{n,X}-EQ_{n,X})^2\le C
B_{f,n}^2.
\end{equation}

\item[(ii)] In addition,
%
\begin{equation}
\label{emaomean} E \Biggl(Q_{n,X}-\sum_{j=1}^{\nu}b_{n,j}f_{X,j}
\Biggr)^2\le C B_{f,n}^2,
\end{equation}
in each of the following three cases.
\begin{enumerate}[(c1)]
\item[(c1)] $d=0$, $g\in
\Lambda_\beta[\Pi]$, $1/2<\beta\le1$;

\item[(c2)] $d\ne0$, $g\in
\Lambda_\beta[\Pi]$,  $1/2<\beta\le1$;

\item[(c3)] $|\dot f_X(u)|\le C u^{-1-2d},$
$ 0< u\le\uppi.$
\end{enumerate}

Moreover, in case \textup{(c1)},
%
\begin{eqnarray}
\label{emeanmean} EQ_{n,X}-\sum_{j=1}^{\nu}b_{n,j}f_j=
\mathrm{o}(B_{f,n}).
\end{eqnarray}
If $b_{n,j}$'s satisfy (\ref{eLF-1}), then
(\ref{emeanmean}) holds also in cases \textup{(c2)} and \textup{(c3)}.
\end{enumerate}
\end{thmm}
\begin{pf} (i) Recall $I_{X,j}=|w_{X, j}|^2. $ By Proposition
\ref{pbound-f},
\begin{eqnarray*}
E(Q_{n,X}-EQ_{n,X})^2= \operatorname{ Var}
\Biggl(\sum_{j=1}^{ \nu} b_{n,j}I_{X,j}
\Biggr) \le C \sum_{j,k=1}^{ \nu}
|b_{n,j}b_{n,k}|\bigl|E[w_{X,j}\overline{w_{X,k}}]\bigr|^2.
\end{eqnarray*}
For $j=k$
bounding $(E|w_{X,j}|^2)^2\le2(E|w_{X,j}|^2-f_{X,j})^2+2f^2_{X,j}$,
and letting
\begin{eqnarray*}
s_{n,1}^\prime &:=&\sum_{j=1}^{ \nu}
b_{n,j}^2\bigl(E|w_{X,j}|^2-f_{X,j}
\bigr)^2, \\ s_{n,2}^\prime&:=&\sum
_{1\le k <j\le\nu} |b_{n,j}b_{n,k}|\bigl| E[w_{X,j}
\overline{w_{X,k}}]\bigr|^2,
\end{eqnarray*}
one obtains
$ 
E(Q_{n,X}-EQ_{n,X})^2\le
C(s_{n,1}^\prime+s_{n,2}^\prime+B_{f,n}^2). $ 
Under the current assumptions, 
by Proposition \ref{propf-mult}(iv), for $1\le
k<j\le\nu$ $(0<\varepsilon<1/2)$,
\begin{eqnarray*}
\bigl |E|w_{X,j}|^2-f_{X,j}\bigr|&\le& Cu_j^{-2d}
\bigl(j^{-1}\log j+\delta_n\bigr),
\\
\bigl|E[w_{X,j}\overline{w_{X,k}}]\bigr|&\le& C \bigl(
\bigl(u_k^{-2d}+u_j^{-2d}\bigr)
j^{-1}\log j 
+ \bigl(u_k^{-2d}
\wedge u_j^{-2d}\bigr) \delta_n\ell(
\varepsilon,j-k) \bigr),
\end{eqnarray*}
where $\delta_n \to0$. Observe that $ s_{n,i}^\prime\le t_{n,i},
i=1,2,$
where $t_{n,1}$ and $t_{n,2}$ are
as in the proof of Lemma~\ref{lemlLip}. 
Therefore, the same argument as
used in proving (\ref{esac-A}) 
implies that $s_{n,1}'+s_{n,2}'$
satisfies the bound (\ref{esac-A}), 
which in turn yields
$s_{n,1}'+s_{n,2}'\le C ( b_{f,n} B_{f,n}+ \varepsilon_{n,G}
B_{f,n}^2 )\le CB_{f,n}^2, $
since $b_{f,n}\le B_{f,n}$.
This completes proof of (\ref{esumsum}).

(ii) By parts (i), (iv) and (iii) of Proposition \ref{propf-mult}, we
respectively obtain
\begin{eqnarray*}
\bigl|E|w_{X,j}|^2-f_{X,j}\bigr|&\le& Cu_j^{-2d}
n^{-\beta},\qquad \mbox{in case (c1)},
\\
&\le& Cu_j^{-2d} \bigl(j^{-1}\log
j+n^{-\beta}\bigr),\qquad \mbox{in case (c2)},
\\
&\le& Cu_j^{-2d} \bigl(j^{-1}\log j\bigr),\qquad
\mbox{in case (c3)}.
\end{eqnarray*}
Let
$D_n:= |EQ_{n,X}-\sum_{j=1}^{\nu}b_{n,j}f_{X,j}  |
= |\sum_{j=1}^{\nu}b_{n,j}(E|w_{X,j}|^2-f_{X,j}) |.
$ 
Under the current assumptions, $f_{X,j}^{-1}\le Cu_j^{2d}$, $0<u\le
\uppi$. Thus, in
case (c1),
%
\begin{eqnarray}
\label{ee1-lk} D_n\le C\sum_{j=1}^{\nu}|b_{n,j}f_{X,j}|n^{-\beta}
\le Cn^{1/2-\beta} \Biggl(\sum_{j=1}^{\nu}(b_{n,j}f_{X,j})^2
\Biggr)^{1/2}=\mathrm{o}(B_{f,n}),
\end{eqnarray}
which proves (\ref{emaomean}) and
(\ref{emeanmean}).

In case (c2),
$ 
D_n\le C\sum_{j=1}^{\nu}|b_{n,j}f_{X,j}|
(j^{-1} \log(n)+n^{-\beta}).
$ 
Arguing 
as for (\ref{earg}), one can show that
$\sum_{j=1}^{\nu}|b_{n,j}f_{X,j}|j^{-1}\log j
=\mathrm{o}(B_{f,n}),$ if (\ref{eLF-1}) holds, and $\sum_{j=1}^{\nu
}|b_{n,j}f_{X,j}|j^{-1}\log j=\mathrm{O}(B_{f,n})$, otherwise,
which
together with (\ref{ee1-lk}) yields (\ref{emaomean}) and~(\ref{emeanmean}).
The proof of (\ref{emaomean}) and~(\ref
{emeanmean}) in case (c3) is
the same as in case (c2). This completes the proof of the theorem.
\end{pf}

\begin{rem}\label{rinterval}
Consider now the sum
%
\begin{eqnarray}
\label{ebnj-small} Q_{n,X}=\sum_{j=1}^{\theta n}b_{n,j}I_{X,j},\qquad
(0<\theta <1/2),
\end{eqnarray}
where summation
is taken over a fraction $\{1, \ldots, \theta n\}$ of the set $\{1,
\ldots,\nu\}$,
and periodograms $I_{X,j}$ used in $Q_{n,X}$ are based on
frequencies $u_j$ from the zero neighborhood $[0,2\uppi\theta]$,
sub-interval of $[0,\uppi]$.
In this case, the smoothness conditions on $f_X$ and $A_X$ are
required only to obtain upper bounds on the covariances
$E[w_{X,j}\overline{w_{X,k}}]$ and $E[w_{X,j}\overline{w_{\zeta,k}}]$
in Proposition \ref{propf-mult}. Therefore, 
in order for
these bounds to be valid at frequencies $u_j\in[0,2\uppi\theta]$ it
suffices to impose smoothness conditions on $f_X$ and $A_X$ on a
slightly larger interval $[0, a]$, $a>2\uppi\theta$, covering
$[0,2\uppi\theta]$. Hence, for the sum $Q_{n,X}$ of (\ref
{ebnj-small}), all of the above results derived in
this section
remain valid if conditions on $f_X$ and $A_X$ are satisfied on some
interval $[0,a]$,
with $a>\Delta$, instead of on $[0,\uppi]$. 
\end{rem}

\begin{rem}\label{rext}
To highlight the main 
method 
of establishing the asymptotic normality of the weighted sums of the
periodograms, we focused mainly on a linear process with an i.i.d. noise~$\{\zeta_j\}$.
Since by the Wold decomposition most of stationary
processes can be written as a linear process with white noise
innovations, 
it is of interest to
extend the above results 
to a linear processes with
martingale-difference innovations. Without assuming that the first
conditional moments of $\zeta_j$ are constant, such extension
requires substantial effort which includes deriving
the general CLT for quadratic forms in martingale-differences
and justification of the Bartlett approximation,
by proving the bound of Proposition \ref{pbound-f}.
Such extension, although non-standard, can be established for a
wide class of martingale difference innovations under tractable
conditions and is currently under our consideration.
\end{rem}

%
\begin{rem}\label{rext1}
In the proof of the asymptotic normality of
the local Whittle
estimator of the memory parameter $d$ in (\ref{effNew}), Robinson
\cite{Rob95N2} established the
CLT
\[
\frac{m^{-1/2}}{g(0)}(S_{n,X} -ES_{n,X})\rightarrow N(0,1),\qquad
S_{n,X}= \sum_{j=1}^m
\nu_{n,j}\frac{I_{X,j}}{|u_j|^{-2d}}= \sum_{j=1}^m
b_{n,j}\frac{I_{X,j}}{f(u_j)}
\]
for special weights $b_{n,j}
=g(u_j)\nu_{n,j}$, where $g(u)\to g(0)>0$, and $\nu_{n,j}:=
\log(j/m) -m^{-1}\sum_{k=1}^m \log(k/m),$ and $m=\mathrm{o}(n)$,
$m\rightarrow\infty$.\vspace*{1pt} Since $ \nu_{n,j}:= \log(j/m)+1+\mathrm{o}(1)$ and
$\sum_{j=1}^m b^2_{n,j}\sim g^2(0)\sum_{j=1}^m \nu^2_{n,j}\sim
g^2(0) m$, they satisfy (\ref{es-con1}) of Theorem \ref{ts-main}
which implies the above CLT. This fact is also apparent upon
examining the Robinson's proof. Additional restrictions on $m$ in
that work were required to show that the bias term
$m^{-1/2}ES_{n,X}$ of the local Whittle estimator is negligible.
\end{rem}

\begin{rem}\label{rtq} 
Here, we provide an example where the weights $b_{n,j}$ in
$S_{n,X}$ do not satisfy
Lindeberg--Feller type condition
(\ref{es-con1}) and the corresponding $S_{n,X}$ does not satisfied
the CLT. Suppose $\{X_j\}$ is a stationary Gaussian zero mean long
memory process, with
$f_X(u)=|u|^{-2d}$, $1/4<d<1/2$.

Let $ \bar X=n^{-1}\sum_{j=1}^n X_j$ and $\hat\gamma(0):=n^{-1}\sum_{j=1}^n (X_j-\bar X)^2$. Recall the identity
\[
2\uppi\sum_{j=1}^nI_X(u_j)=
\sum_{j=1}^n X_j^2=
\sum_{j=1}^n(X_j-\bar
X)^2+ n \bar X^2=n\hat\gamma(0)+ n \bar X^2.
\]
%

Suppose $n$ is even and $\nu=n/2 -1$. Since
$I_X(u_j)=I_X(u_{n-j}), 1\le j\le n$, and $2\uppi I_{X}(u_0)=n\bar
X^2$, we obtain $2\uppi\sum_{j=1}^nI_X(u_j)=4\uppi
\sum_{j=1}^\nu I_X(u_j)+$ $ 2\uppi\{I_X(u_0)$ $+I_X(u_{n/2})\}$, and
\[
 4\uppi\sum_{j=1}^\nu
I_X(u_j)=n\hat\gamma(0) - 2\uppi I_X(u_{n/2}).
\]
%
Now, let $b_{n,j}:=n^{-2d}
4\uppi f_X(u_j)=4\uppi(2\uppi j)^{-2d}$. Then
\[
S_{n,X}= \sum_{j=1}^\nu
b_{n,j} \frac{I_X(u_j)}{f_X(u_j)}= n^{-2d}4\uppi\sum
_{j=1}^\nu I_X(u_j).
\]
%

By Hosking (\cite{Hos96}, Theorem 4),
under the assumed set up here, $n^{1-2d}(\hat\gamma(0)-E\hat\gamma
(0))\rightarrow_D Y$, where
$Y$ is a non-Gaussian r.v. Arguing as in the proof of Theorem
\ref{thmbias-t}(i), one can verify that
$\operatorname {Var}(I_X(u_{n/2}))=\mathrm{O}(1)$. Hence, $ S_{n,X}-ES_{n,X}
\rightarrow_D Y $ does not satisfy the CLT. It remains to show
that for $d>1/4$, the weights $b_{n,j}$ do not satisfy
(\ref{es-con1}):
\[
\frac{\max_{j=1,\ldots,
\nu}|b_{n,j}|^2}{\sum_{j=1}^{\nu}b^2_{n,j}} 
=\frac{\max_{j=1,\ldots, \nu}|j|^{-4d}}{\sum_{j=1}^{\nu}j^{-4d}} \rightarrow
\frac{1}{\sum_{j=1}^{\infty}j^{-4d}}>0.
\]

Moreover, Theorem \ref{ts-main} does not provide the asymptotic
of $\operatorname{ Var}(S_{n,X})$ and approximations of Lemma \ref{lemlap}
break down. To see that, we now have
$b_n=2(2\uppi)^{1-2d}$,
$B_n^2=q_n^2=\sum_{j=1}^nb_{n,j}^2\rightarrow4(2\uppi
)^{2-4d}\sum_{j=1}^\infty j^{-4d}$. 
Since the error of approximation 
$E(S_{n,X}-S_{n,\zeta})^2\le
C\log^{4}(n)$ in (\ref{esacon}) is no more negligible compared to
$B_n^2$, the claim that $\operatorname{ Var}(S_{n,X})\sim q_n^2$ of
Theorem
\ref{ts-main} does not hold.
On the other hand, by Theorem 3 of Hosking \cite{Hos96},
$\operatorname{ Var}(n^{1-2d} \hat\gamma(0))\rightarrow
C>0$,
so that
$\operatorname{ Var}(S_{n,X})\rightarrow
C$.
\end{rem}


\begin{ex}\label {exARMA}Consider the stationary
$\operatorname{ARFIMA}(p,d,q)$ model
\begin{eqnarray*}
\phi(B)X_j=(1-B)^{-d}\theta(B)\zeta_j,\qquad  j\in
\Z,\qquad \{\zeta_j\}\sim \operatorname{IID} \bigl(0,\sigma^2_\zeta
\bigr),\qquad |d|<1/2.
\end{eqnarray*}
Hosking
\cite{Hos81} has shown that the spectral density $f_X$ of this model
satisfies (\ref{effNew}).
We shall show it also satisfies 
(\ref{eaaa}).
Let $h(u)= (1-\mathrm{e}^{-\i u })^{-d}$ and $A_Y(u)=\theta(\mathrm{e}^{-\i u
})/\phi(\mathrm{e}^{-\i u })$.
The transfer function $A_X$ can be written as
%
\begin{eqnarray}
\label{eargum} 
A_X(u)=h(u)A_Y(u),\qquad
f_X(u)=\bigl|A_X(u)\bigr|^2.
\end{eqnarray}
%
Now observe that
$h$ is differentiable and satisfies $
|h(u)|\le C|u|^{-2d},$
$|\dot{h}(u)|\le C|u|^{-1-2d},$ for all $u\in[0,\uppi],$ and
$|h(u)|\sim|u|^{-2d},$ as $ u\to0.$
Thus, for all $0<|u|<\uppi$,
\begin{eqnarray*}
\bigl|\dot{A}_X(u)\bigr|&\le& C\bigl(\bigl|\dot{ h}(u)\bigr|\bigl|A_Y(u)\bigr|+\bigl|h(u)
\bigl|\bigl|\dot{A}_Y(u)\bigr|\bigr) \le C\bigl|1-\mathrm{e}^{-\i u}\bigr|^{-d-1} \le
C|u|^{-d-1},
\end{eqnarray*}
and hence $A_X$ satisfies (\ref{eaaa}). Note also that
$A_X=hA_Y$ is naturally factored
into a differentiable component $h$ and continuous
component $A_Y$ as required in Theorem
\ref{tcontinuous1}. Thus, Theorems
\ref{ts-main}, \ref{tcontinuous}--\ref{thmbias-t} are applicable.
\end{ex}

\begin{ex} Now consider a more general process $\{X_j\}$,
\begin{eqnarray*}
X_j=(1-B)^{-d}Y_j, \qquad j\in\Z,\ |d|<1/2,
\end{eqnarray*}
%
where $ Y_j= \sum_{k=0}^\infty b_k\zeta_{j-k}, $ $ \{\zeta_j\}\sim \operatorname{IID}
(0,1), \sum_{k=0}^\infty|b_k|<\infty, $ is a short memory
process.
Because, letting $A_Y(u)=\sum_{k=0}^\infty b_ke^{-\i uk}$,
$f_X$ and $A_X$
are the same as in (\ref{eargum}), the same argument as used
in Example \ref{exARMA} shows that $f_X$ satisfies (\ref{effNew})
with parameter $|d|<1/2$. Although $A_X$ may not satisfy
(\ref{eaaa}), because $A_Y$ is only continuous, but $A_X$ is
factored as required in Theorems \ref{tcontinuous1} and
\ref{thmbias-t}. Hence, these two theorems are applicable.
\end{ex}
%

\section*{Acknowledgements}
Research supported in part by ESRC
Grant RES062230790 and the NSF Grant DMS-07-04130. Authors would like
to thank four reviewers of the paper for their constructive
comments.


%


%
%

\printhistory

\end{document}